\documentclass[10pt]{iopart}
\long\def\mod#1{|#1|}
\usepackage{amssymb,amstext}
\usepackage[dvips]{graphicx}

\def\Rset{{\cal R}}
\def\zmax{{\ell(\tau)}}
\def\zmaxa[#1]{{\ell(#1)}}

\def\Tim{{T^*}}
\newcommand\R{{\mathbb R}}
\newcommand\N{{\mathbb N}}
\newcommand{\eop}{\hspace*{3mm}\hbox{}\nobreak\hfill\hbox{\vrule height7pt width7pt}\par\medbreak}
\def\qed{\relax\ifmmode\hskip2em \Box\else\unskip\nobreak\hskip1em $\Box$\fi}

\newtheorem{teo}{Theorem}[section]

\newtheorem{coro}[teo]{Corollary}

\newtheorem{theorem}[teo]{Theorem}
\newtheorem{lemma}[teo]{Lemma}
\newcommand{\dis}{\displaystyle}
\begin{document}
\title[Self-similar blow-up]{Self-similar  blow-up for a diffusion-attraction problem}
\author{Ignacio A. GUERRA\dag\footnote[3]{Supported by FONDECYT 3040059}, Mark A. PELETIER\ddag}

\address{\dag Centre for Mathematical Modeling, 
Universidad de Chile, FCFM \\
Casilla 170 Correo 3 Santiago, Chile}
\address{\ddag Department of Mathematics and Computer Science, TU Eindhoven, \\
P.O Box 513 5600 MB Eindhoven, The Netherlands}

\eads{\dag\mailto{iguerra@dim.uchile.cl}}
\eads{\ddag\mailto{Mark.Peletier@cwi.nl}}

\begin{abstract}
In this paper we consider a system of equations that describes
a class of mass-conserving aggregation phenomena, including
gravitational collapse and bacterial chemotaxis.
In spatial dimensions strictly larger than two, and under the assumptions
of radial symmetry, it is known that this system has at
least two stable mechanisms of singularity formation (see
e.g.\ M.~P.\ Brenner et al.\ 1999, {\it Nonlinearity} {\bf 12},
1071-1098); one type is self-similar, and may be viewed as a
trade-off between diffusion and attraction, while in the other
type the attraction prevails over the diffusion and a non-self-similar
shock wave results. Our main
result identifies a class of initial data for which the
blow-up behaviour is of the former, self-similar type.
The blow-up profile is characterized as belonging to a subset of
stationary solutions of the associated ordinary differential
equation.
\end{abstract}
\submitto{\NL}
\section{Introduction}
We consider
the parabolic-elliptic system
\begin{eqnarray}
\dis n_t={\rm div}\{\Theta\nabla n+n\nabla \phi\} \quad
&\mbox{in}\quad \Omega\times \R^+,
\label{P1.s}
\\
%
%
\dis \Delta \phi=
n\qquad\qquad\qquad\qquad\quad\,\,&\mbox{in}\quad \Omega\times
\R^+,
\label{P2.s}
\\
\dis 0=(\Theta\nabla n+n\nabla \phi)\cdot \vec \nu\quad
\,\,&\mbox{on}\quad\partial \Omega\times \R^+,
\label{P5.s}
\\
\dis \phi=0 \qquad \qquad \qquad\qquad\quad\:\:
&\mbox{on}\quad\partial \Omega\times \R^+,
\label{P4.s}
%
%
%
\\
 \dis
n(x,0)=n_0(x)\qquad\qquad\qquad\quad \:\:\: &\mbox{in}\quad\Omega,
\label{P6.s}
\end{eqnarray}
where $\Omega=B_1(0)=\{ x\in \R^d\colon |x|\leq 1\},$  $d>2,$ and
$\vec\nu$ is the outer normal vector from the boundary
$\partial\Omega.$ Here $\Theta>0$ is a constant parameter. The
initial condition $n_0$ is chosen in $L^2(\Omega)$, radially
symmetric, and such that
\begin{eqnarray}\label{mass1int.s}
\int_\Omega n_0 \,dx=1,\quad\mbox{and}\quad n_0(x)\geq 0\quad
\mbox{in}\quad\Omega.
\end{eqnarray}
Equations (\ref{P1.s})--(\ref{mass1int.s}) define a problem for
the unknown mass density $n$ and potential $\phi$. 
Mass is conserved  by the no-flux condition (\ref{P5.s}), and 
therefore~(\ref{mass1int.s}) implies
\begin{eqnarray}\label{unitarymass}
\int\limits_{\Omega}n(x,t)\,dx=\int\limits_{\Omega}n_0(x)\,dx=1.
\end{eqnarray}

Problem (\ref{P1.s})--(\ref{mass1int.s}) is  a model for the
evolution of a cluster of particles under gravitational
interaction and Brownian motion (see \cite{bkn} and the references
therein). Here $n$ represents the mass density, $\phi$ the gravitational
potential, and $\Theta$ a rescaled temperature characterizing the Brownian
motion. 
%
This model also appears
in the study of evolution of polytropic stars,  by considering the
evolution of self-interacting clusters of particles under
frictional and fluctuating forces \cite{w2}.
Finally, problem (\ref{P1.s})--(\ref{mass1int.s}) also arises in the study
of the motion of bacteria by chemotaxis as a simplification (see
\cite{jl}) of the Keller-Segel model \cite{ks,v,bb,blb}. Here the variables
$n$ and $\phi$ represent the density of bacteria and the
concentration of the chemo-attractant.

We view the problem~(\ref{P1.s})--(\ref{mass1int.s})
as an evolution equation in $n$, since
by equations~(\ref{P2.s}-\ref{P5.s}) the
function $\phi$ is readily recovered from the solution $n$.
It is known~\cite{bn2} that problem
(\ref{P1.s})--(\ref{mass1int.s}) has a unique local solution
if $n_0\in L^2(\Omega)$, which satisfies
$n\in L^\infty\big(\Omega\times (\epsilon,\tilde T)\big)$ for some
$\tilde T>0$ and for every $\epsilon > 0$.
We restrict ourselves to
the analysis of radially symmetric solutions and write $n(r,t):=n(x,t)$ with
$r=|x|\in [0,1].$

Since we are interested in the question when and how~(\ref{P1.s})--(\ref{mass1int.s})
generates singularities, we define:
\[
\fl\qquad\Tim=\sup\{\, \tau>0\mid \mbox{Problem
(\ref{P1.s})--(\ref{mass1int.s}) has a solution}\,n \in
L^\infty(\Omega\times(\epsilon,\tau])\,\}.
\]
If $\Tim<\infty$, then we say that blow-up occurs for
(\ref{P1.s})--(\ref{mass1int.s}), in which case
\begin{eqnarray}\label{blowlimn}
\lim\limits_{t\to \Tim}\sup\limits_{[0,1]} n(r,t)=\infty.
\end{eqnarray}
Various sufficient conditions for blow-up are known~\cite{b1,bhn,bn1,bn2}.

\begin{figure}[htb]
\begin{center}
\includegraphics[angle=0,width=2.5in,height=1.8in]{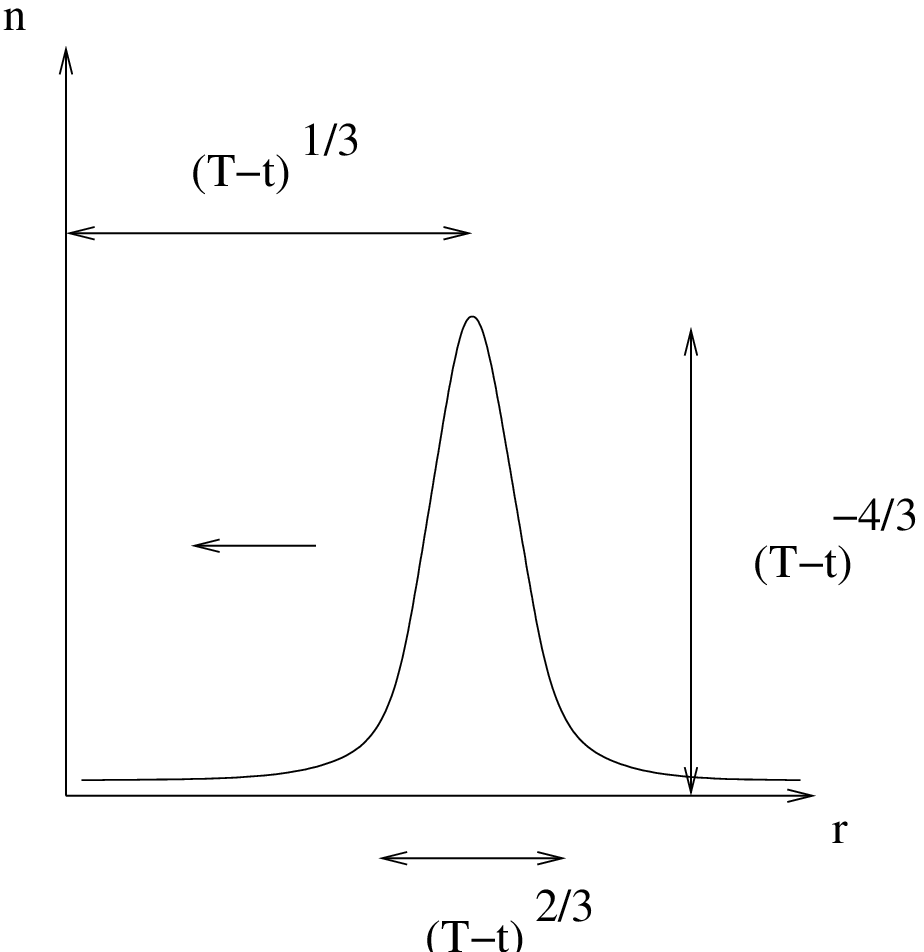}
\includegraphics[angle=0,width=2.5in,height=1.8in]{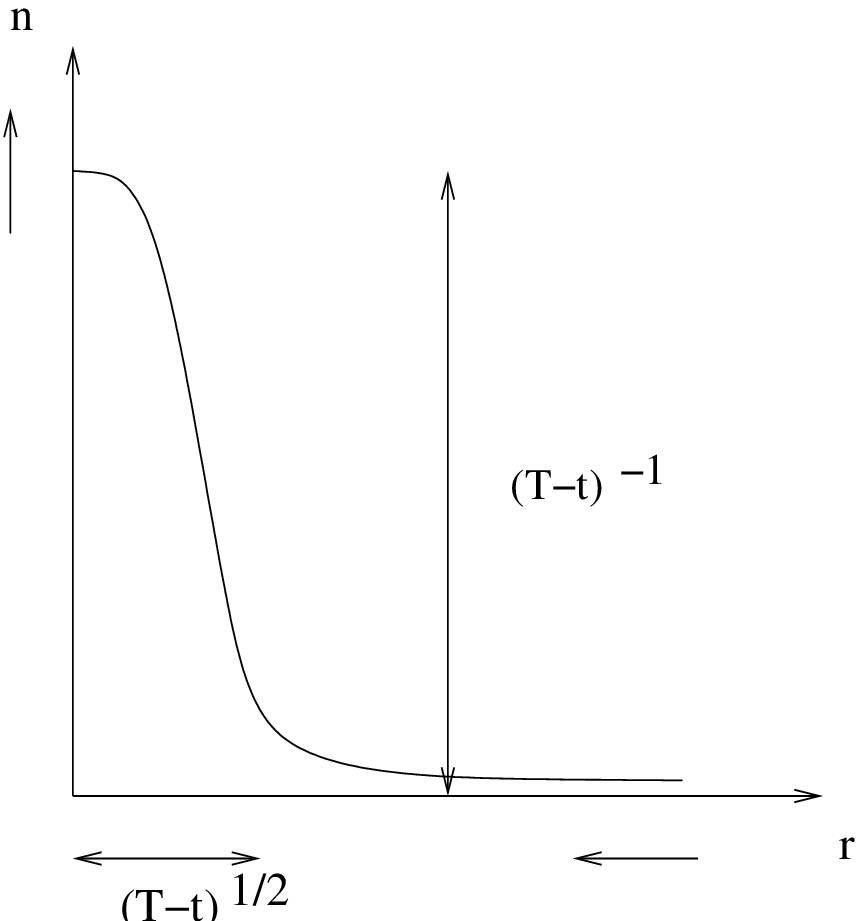}
\end{center}
\caption{The profile $n(r,t)$ for blow-up with(left) and
without(right) con\-cen\-tra\-tion of mass, with $T=\Tim.$
\label{fiblow}}
\end{figure}

For $d=3,$ Herrero et al.~\cite{hmv1,hmv2} were the first to study the
behaviour of the solution close to blow-up, using matched asymptotic expansions.
Later Brenner et al.~\cite{bcksv} studied the problem for
$2<d<10.$ They used a numerical approach to describe solutions and
proved existence and linear stability of similarity
profiles. Note however that no proof of convergence or
characterization of blow-up in terms of initial data were given in
these references.
The principal types of blow-up described in \cite{hmv1,hmv2,bcksv}
are:
\begin{enumerate}
\item[(a)] A solution $n(r,t)$  consists of an imploding
smoothed shock wave which moves towards the origin.  As $t\to
\Tim,$ the bulk of such a wave is concentrated at distances
$O((\Tim-t)^{1/d})$ from the origin, has a width
$O((\Tim-t)^{(d-1)/d}),$ and at its peak it reaches a height of
order $O((\Tim-t)^{-2(d-1)/d}).$ This type of blow-up has the
property of concentration of mass at the origin at the blow-up
time, i.e.
\begin{eqnarray}\label{concmassint}
\lim\limits_{r\to 0}\left[\lim\limits_{t\to \Tim}\int\limits_{0}^r
n(y,t)y^{d-1}\,dy\right]=C>0.
\end{eqnarray}
This situation is depicted in Figure \ref{fiblow} (left). 
\item[(b)] A solution $n(r,t)$ has a self-similar blow-up of the
form
\begin{equation}
\label{sim:sss}
(\Tim-t){\dis n}\left(\eta\sqrt{(\Tim-t)\Theta},t\right)\sim
\Psi(\eta)
\quad \mbox{as}\quad t\to \Tim.
\end{equation}
Note that this implies that $n$ satisfies (\ref{concmassint}) with
$C=0.$ Therefore no concentration of mass at the origin occurs at
the blow-up time. This blow-up behaviour is depicted in Figure
\ref{fiblow} (right).
\end{enumerate}

\bigbreak

The results of this paper are two-fold. First, we demonstrate rigorously that 
the self-similar blow-up structure~(\ref{sim:sss}) is an attractor for
the system~(\ref{P1.s})--(\ref{mass1int.s});
secondly, we identify an explicit class of initial data that converges to 
a self-similar solution of this type. Let us elaborate on this.

Let
$n_0=n_0(r)$ be such that
\begin{eqnarray}\label{inicon0.n}
\fl \chi_d r^d n_0(r)\leq\|n_0\|_{L^1(B_r(0))}\:\:\mbox{for}
\:\:r\in (0,1),\\
\fl\Theta(n_0)_r+n_0(\phi_0)_r\geq 0,\quad(r^d(\phi_0)_r)_r=r^d
n_0\:\:\mbox{in}\:\:(0,1),\quad \mbox{and}\:\:\phi_0(1)=0,
\label{btpositive.n}
\end{eqnarray}
where $\chi_d$ is the measure of the unit ball in $\R^d$.
Suppose also that $\Theta\leq 1/(4d\chi_d)$, implying
that the solution $n=n(r,t)$ of (\ref{P1.s})--(\ref{mass1int.s}) blows up
at finite time $\Tim>0$ and at the point $r=0$~\cite{bhn}. Finally, assume that
the two functions 
\begin{eqnarray}\label{inter.n}
\fl\|n_0\|_{L^1(B_r(0))}\quad \mbox{and}\quad \frac{4 \Theta
r^d}{2(d-2)\Theta\Tim+r^2}\:\:\mbox{intersect exactly once in}\:\:
[0,1].
\end{eqnarray}
Our main result (Theorem~\ref{prin.s2}) shows that if
(\ref{inicon0.n}), (\ref{btpositive.n}), and (\ref{inter.n}) hold,
then $n$ satisfies
\[
n(0,t)\leq \frac{2d}{(d-2)}(\Tim-t)^{-1}\quad \mbox{for}\quad t\in
(0,\Tim),
\]
and moreover has a structure near blow up given by
\[
n_*(r,t)=(\Tim-t)^{-1}\Psi\left(\frac r{
\sqrt{\Theta(\Tim-t)}}\right),
\]
where the function $\Psi$ is one of a class of solutions of a
steady-state problem; a class that includes the functions
\[
\Psi_1(\eta):=
(d-2)\frac{(2d+\eta^2)}{(d-2+\frac{1}{2}\eta^2)^2}\quad
\mbox{and}\quad \Psi^*(\eta):= 1\quad \mbox{for}\quad\eta>0.
\]
In particular the initial state  $n_0\equiv 1/\chi_d$ and
$\Theta\leq 1/(4d\chi_d)$
satisfies the conditions above (Corollary~\ref{prin.s3}). If we
relax assumption~(\ref{inter.n}) but assume instead that $n$
satisfies the growth condition
\[
n(0,t)\leq M(\Tim-t)^{-1}\quad \mbox{for}\quad t\in (0,\Tim),
\]
for some constant $M>0,$ then $n$ has the same structure of blow-up 
given above (Theorem \ref{prin.s}). The hypotheses on the
initial data (\ref{inicon0.n}), (\ref{btpositive.n}), and
(\ref{inter.n}) are more natural in the context of a transformed problem we
introduce in the next section. Note however that $(n_0)_r\leq 0$
in $[0,1]$ implies assumption (\ref{inicon0.n}).

This paper is organized as follows. In section \ref{mresults}, we
put the problem in terms of a new variable, thus transforming the
system (\ref{P1.s})--(\ref{mass1int.s}) into a single PDE, and then
state our results in terms of this new formulation. In section
\ref{discussion}, we discuss some non-self-similar blow-up patterns
related to case (a). Sections \ref{prel}, \ref{convtheo}, and \ref{intersectioncomparison}
provide the tools for the proofs of Theorems \ref{prin.s2} and~\ref{prin.s}, and the
arguments are 
wrapped up in Section~\ref{proofresults}.
A rather technical dervation of a Lyapunov function is placed in Appendix~\ref{sec:appendix},
and in Appendix~\ref{stability} we derive some linear stability results.

\section{Precise statements of main results}\label{mresults}
For radial solutions, the average density function $b(r,t)$~\cite{bcksv} is
defined  by
\begin{eqnarray}\label{newdens}
b(r,t):=\frac{d\,\chi_d}{r^d}\int\limits_{0}^r n(y,t)y^{d-1}\,dy,
\end{eqnarray}
This variable turns out to be convenient in the analysis of this
system.
Note that it has the same scale invariance  as
$n(r,t)$, but that solutions are smoother when expressed in terms
of $b$. For example, if for some fixed $t>0$ the density
$n(r,t)$ is a delta function at the origin with unit mass, then
$b(r,t)=r^{-d}.$
Let $D=(0,1)$ and set $D_T=D\times (0,T)$ for some $T>0.$
Transformation (\ref{newdens}) puts system
(\ref{P1.s})--(\ref{mass1int.s}) in the form
\begin{eqnarray}\label{bequation}
b_t=
\chi_d\Theta\left(b_{rr}+\frac{d+1}{r}b_r\right)+\frac
1{d}rbb_r+b^2\qquad&\mbox{in}\quad D_T \\
\label{bequation1}
b_r(0,t)=0,\quad b(1,t)=1,& \mbox{for}\quad t\in [0,T), \\
\label{bequation2} b(0,r)=b_0(r) & \mbox{for}\quad r\in D.
\end{eqnarray}
Here we have redefined $t:=\frac 1{\chi_d}t.$ Regarding the
initial condition, we assume
\begin{eqnarray}\label{inic}
b_0\in C^2(\overline D),\quad\mbox{and}\quad \frac r{d}(b_0)_r+b_0\geq
0\quad \mbox{for}\quad r\in D,
\end{eqnarray}
where the second condition is equivalent to $n_0\geq 0$ in $D.$
Note that the conservation of the mass (\ref{unitarymass}) is
represented by $b(1,t)=1$ for $t\in [0,T).$
As was done for problem (\ref{P1.s})--(\ref{mass1int.s})   we
define $T>0$ to be the maximal time of existence for the average
density $b(r,t).$ If $\Tim<\infty$ in (\ref{blowlimn}),
then
\[
\lim\limits_{t\to T}\sup\limits_{[0,1]} b(r,t)=\infty,
\]
where $T=\Tim/\chi_d.$ Using (\ref{newdens}), we deduce $
b(r,t)\leq 1/r^d$ for $r\in \overline{D},\:t>0;$ this implies
single point blow-up for $b(r,t)$ at the point $r=0$.
To characterize  the asymptotic behaviour near blow-up of the
solution $b(r,t)$ of problem (\ref{bequation})--(\ref{inic}),
we study the solutions of the associated
boundary-value problem
\begin{eqnarray}\label{staW.n}
\cases{ \varphi_{\eta\eta}+\frac {d+1}{\eta}\varphi_\eta+\frac
1{d}\eta \varphi \varphi_\eta- \frac 1{2}\eta
\varphi_\eta+\varphi^2-\varphi=0, \quad\mbox{for}\quad  \eta>0,
\cr \varphi(0)\geq 1\quad \varphi_\eta(0)=0. }
\end{eqnarray}
If  $b$ is a solution  of (\ref{bequation})--(\ref{inic}) which
blows up at time $T>0$ and at the point $r=0$, then we will show
that it has the asymptotic form given by
\[
b_*(r,t)=(T-t)^{-1}\varphi\left(\frac r{
\sqrt{\chi_d\Theta(T-t)}}\right).
\]
Equation (\ref{staW.n}) has multiple solutions for $2<d<10$~\cite{hmv2,bcksv}.
We classify them by counting the number of
times they cross the singular solution
$\varphi_S(\eta):=2d/\eta^2.$ For that purpose, we introduce the
set
\[
{\cal S}_k=\{  \varphi \colon\;\mbox{$\varphi$ is a solution of
(\ref{staW.n}) that has $k$ intersections with $\varphi_S$}\}.
\]
We shall see that ${\cal S}_1$ is the relevant subset of solutions
of (\ref{staW.n}) for the characterization of the type of blow-up considered in
this paper. Numerical
evidence \cite{bcksv} suggests that ${\cal S}_1$ contains only two
elements: \begin{equation}\label{phis1} \varphi^*(\eta)= 1\quad \mbox{and}\quad
\varphi_1(\eta):=\frac{2d}{(d-2+\frac{\eta^2}{2})}\quad\mbox{for}\quad
\eta\geq 0.
\end{equation}
For the initial condition, we assume
\begin{eqnarray}\label{inicon0}
(b_0)_r\leq 0\quad \mbox{for}\quad r\in D,
\end{eqnarray}
and
\begin{eqnarray}\label{btpositive}
\chi_d\Theta\left((b_0)_{rr}+\frac{d+1}{r}(b_0)_r\right)+\frac
1{d}r b_0\,(b_0)_r+b_0^2\geq 0\quad \mbox{for}\quad r\in D.
\end{eqnarray}
We will show that this implies $b_r\leq 0$ in $D_T$ and $b_t\geq
0$ in $D_T.$ In terms of $n_0$ assumption~(\ref{inicon0}) becomes
(\ref{inicon0.n}) and assumption (\ref{btpositive}) becomes
(\ref{btpositive.n}).
\begin{teo}\label{prin.s2}
Let $d>2$ and $b_0$ satisfy (\ref{inicon0}) and
(\ref{btpositive}). Let $b(r,t)$ be the corresponding solution of
problem (\ref{bequation})--(\ref{inic}) that blows up at $r=0$ and
at $t=T.$ If
\begin{eqnarray}\label{corbB.1}
\fl\quad\Theta\leq \Theta_1 := 1/(4d\chi_d) \quad\mbox{and}\quad
\quad b_0(r)\ \mbox{intersects} \
T^{-1}\varphi_1(r/\sqrt{\chi_d\Theta T})\ \mbox{once}
\end{eqnarray}
then
\begin{eqnarray}\label{corobB.2}
b(0,t)\leq  M_1 (T-t)^{-1}\quad \mbox{for}\quad
t\in(0,T)
\end{eqnarray}
with $M_1:= 2d/(d-2)$. Moreover,  $T<M_1/b_0(0)$, and there exists $\varphi\in S_1$ such
that
\begin{eqnarray}\label{selfsimconv.n}
\lim\limits_{t\to
T}(T-t)b\left(\eta\sqrt{\chi_d\Theta(T-t)}\right)=\varphi(\eta)
\end{eqnarray}
uniformly on compact sets $|\eta|\leq C$ for every $C>0.$
\end{teo}
We remark that there exists a family of $b_0$ satisfying the
conditions (\ref{inic}), (\ref{inicon0}), and  (\ref{btpositive}),
given by $b_0(r)=K_1+K_2/(r^d+K_3)$ with positive constants $K_i$ that
satisfy $K_1+K_2/(1+K_3)=1$ and $\Theta<K_2/2 d^2\chi_d.$
Conditions (\ref{inic}), (\ref{inicon0}), and  (\ref{btpositive})
are also satisfied for $b_0\equiv1$.
Note that condition (\ref{corbB.1}) of Theorem \ref{prin.s2} can be
generalized by changing $\varphi_1$ for other
solution $\varphi$ of (\ref{staW.n}). Since these solutions are only
known numerically, the counterpart of $M_1$ and $\Theta_1$ cannot be given explicitly.
The next corollary applies this result to
$b_0\equiv 1.$
\begin{coro}\label{prin.s3}
Let $d>2,$ $b_0\equiv 1$, and $\Theta< \Theta_1.$ Then $b(r,t)$,
the corresponding solution of problem
(\ref{bequation})--(\ref{bequation2}), blows up at $r=0$ and at
some time $t=T<M_1$; moreover (\ref{corobB.2}) holds  and there
exists $\varphi\in S_1$ satisfying~(\ref{selfsimconv.n}).
\end{coro}
Numerical simulations \cite{bcksv} suggest that for an open set of
initial data the convergence in (\ref{selfsimconv.n}) holds for
$\varphi=\varphi_1.$ This self-similar behaviour may be seen roughly in
Figure~\ref{fiblow} (right), by imagining  $n(r,t)$ replaced by $b(r,t)$
(since $n$ and $b$ scale similarly).
In \ref{stability} we show that $\varphi_1$ is linearly stable
(using the result in \cite{bcksv}) and also that $\varphi^*$ is
linearly unstable.

For more general initial data we have the following result.
\begin{teo}\label{prin.s}
Let $d>2$ and let $b_0$ satisfy (\ref{inicon0}) and
(\ref{btpositive}). Assume that $b(r,t),$
 the corresponding solution of problem
(\ref{bequation})--(\ref{inic}),  blows up at $r=0$ and at $t=T.$
If $b$ satisfies the growth condition
\begin{eqnarray}\label{corobB}
b(0,t)\leq M (T-t)^{-1}\quad \mbox{for}\quad t\in(0,T)
\end{eqnarray}
with $M>0,$ then there exists $\varphi\in S_1$ such that the
convergence (\ref{selfsimconv.n}) holds.
\end{teo}

We now briefly discuss the structure of the proofs of these theorems. 
Following the scale invariance, we set
\[
\fl \tau=\log\left(\frac{T}{T-t}\right), \quad 
\eta=\frac{r}{(\chi_d\Theta(T-t))^{1/2}},
\quad \text{and}\quad B(\eta,\tau)=(T-t)b(r,t).
\]
The rectangle $D_T$ transforms into
\[
\fl\qquad\Pi=\{(\eta,\tau)\mid \tau>0,
\,0<\eta<\zmax\,\}\quad\mbox{where}\quad \zmax:=(\chi_d\Theta
T)^{-1/2}e^{\tau/2}.
\]
The initial-boundary problem (\ref{bequation})--(\ref{inic}) now
becomes
\begin{eqnarray}
\label{equaB1}  \fl\qquad B_\tau+B+\frac 1{2}\eta
B_\eta=B_{\eta\eta}+\frac
{d+1}{\eta} B_\eta+\frac 1{d}\eta BB_\eta+B^2\quad & \mbox{in}\:\Pi, \\
\fl\qquad B_\eta(0,\tau)=0,\quad
B\left(\zmax,\tau\right)=e^{-\tau}T&
\mbox{for}\: \tau\in\R^+, \\
\label{equaB3}
 \fl\qquad B(\eta,0)=B_0(\eta):=T b_0\left(\eta (\chi_d\Theta
T)^{1/2}\right) &\mbox{for}\:\eta\in \Pi(0),
\end{eqnarray}
where $\Pi(0)=(0,\zmaxa[0]).$
Note that a solution of (\ref{staW.n}) is a time-independent
solution of (\ref{equaB1})--(\ref{equaB3}). Therefore the study of
the  blow-up behaviour of $b(r,t)$ is reduced to the analysis of the
large time behaviour of solutions $B(\eta,\tau)$ of
(\ref{equaB1})--(\ref{equaB3}), and in particular stabilization
towards solutions $\varphi$ of (\ref{staW.n}).
The proof of Theorem \ref{prin.s} consists of two parts. In
Section  \ref{convtheo}, we first prove that
$\omega\subset {\cal S}_1,$ where
\begin{eqnarray}\nonumber
\fl\omega=\{\phi\in L^\infty(\R^+) : 
\exists \tau_j\to \infty \quad\mbox{such that}\quad \\
\fl \qquad B(\cdot,\tau_j)\to \phi(\cdot)\quad \mbox{as}\quad
\tau_j\to \infty \quad \mbox{uniformly on compact subsets of}\quad
\R^+\}
\label{omegalimitset}
\end{eqnarray}
is the  $\omega$-limit set we introduce for
(\ref{equaB1})--(\ref{equaB3}). The proof
uses the observation that
equation (\ref{equaB1}), without the convection
term $\frac 1{d}\eta BB_\eta,$ is the backward self-similar
equation for the parabolic semilinear
equation
\begin{eqnarray}\label{semi-linear}
\bar b_t=\Delta_N \bar b +\bar b^2,
\end{eqnarray}
where $\Delta_N$ denotes the Laplacian in $\R^N$ and
$N=d+2$~\cite{gk0,gk}. We use the methods for the analysis of this
self-similar equation to prove Theorem \ref{prin.s}. However, due
to the presence of the convection term, a different Lyapunov
functional is necessary. This functional is constructed using the
method of Zelenyak~\cite{zel}, which yields a Lyapunov functional
in implicit form. In section \ref{intersectioncomparison}, we use
intersection comparison arguments based on the ideas of
Matano~\cite{mat} to prove that the $\omega$-limit set
(\ref{omegalimitset}) is a singleton. With a result on
intersection with $\varphi_S$ this completes the proof of Theorem
\ref{prin.s}.

Note that Theorem \ref{prin.s} is similar to a result for the
supercritical case $(N>6)$ for equation (\ref{semi-linear}), where
two different kinds of self-similar blow-up behaviour may coexist
\cite{m}.

Finally to obtain Theorem \ref{prin.s2} and Corollary
\ref{prin.s3}, we use Theorem \ref{prin.s} and comparison ideas
from Samarskii et al. \cite[Chapter IV]{sgkm}.

\section{Discussion on non self-similar blow-up patterns}\label{discussion}

In this section we discuss a family of blow-up patterns which
appears when we refine the asymptotic expansion for the profile
$\varphi=\varphi^*\equiv 1.$  This situation is closely related to
the blow-up behaviour of  (\ref{semi-linear}) with $N<6.$ If a
solution $\bar b$ of (\ref{semi-linear}) with $N<6$
blows up at $x=0$ and $t=T$, then
\[
\lim\limits_{t\to T}(T-t)\bar b(\eta\sqrt{T-t},t)=1
\]
uniformly on compact sets $|\eta|<C$ for arbitrary $C>0$~\cite{gk0,gk}.
Moreover it have been shown (see for instance  \cite{m1,vel}) that
a refined description of blow-up gives the existence of two
possible types of behaviour: either
\begin{eqnarray}\label{semiphi1}
\lim\limits_{t\to T}(T-t)\bar
b\left(\eta\sqrt{(T-t)|\log(T-t)|},t\right)= \bar \varphi_1(\eta)
\end{eqnarray}
uniformly on compact sets $|\eta|<C,$ with $C>0$ arbitrary; or
\begin{eqnarray}\label{semiphim}
\lim\limits_{t\to T}(T-t)\bar
b\left(\eta(T-t)^{1/2m},t\right)=\bar \varphi_m(\eta) \quad
\mbox{for some}\quad  m\geq 2,
\end{eqnarray}
uniformly on compact sets $|\eta|<C,$ with $C>0$ arbitrary. Here
the family $\{\bar\varphi_i\}_{i\geq 1}$ is known explicitly.
For problem (\ref{bequation})--(\ref{inic}), it was shown~\cite{hmv2}
for $d=3$ that there exists a refined asymptotics
for $\varphi^*\equiv 1$. Extending the argument to all $d>2,$
these asymptotics suggest a convergence given by either
\begin{eqnarray}\label{chemphi1}
\lim\limits_{t\to T} (T-t)
b\left(\eta\sqrt{(T-t)|\log(T-t)|^{(d-2)/d}},t\right)=\tilde\varphi_1(\eta)
\end{eqnarray}
or
\begin{eqnarray}\label{chephim}
\lim\limits_{t\to T} (T-t) b\left(\eta(T-t)^{\frac 1{d}+
\frac{d-2}{2(m+d-1)}},t\right)= \tilde\varphi_m(\eta)
\end{eqnarray}
for some $m\geq 2.$
An implicit  formula for the family $\{\tilde\varphi_m\}_{m\geq
1}$ is given in \cite[equation~(43)]{bcksv}.
The type of convergence in $\eta$ towards these profiles is an
open problem.
In (\ref{chephim}), we can take formally the limit $m\to \infty$
and find a non-trivial scaling,
\begin{eqnarray}\label{phity}
\lim\limits_{t\to T}(T-t)
b(\eta(T-t)^{1/d},t)=\tilde\varphi_\infty(\eta).
\end{eqnarray}
Note that this limit cannot be taken for the semilinear equation
where (\ref{semiphim}) holds. The convergence (\ref{phity})
represents the convection-dominant behaviour of
(\ref{bequation})--(\ref{inic}), which in terms of the density
$n=n(r,t)$ describes an imploding wave moving towards the origin,
as shown in Figure \ref{fiblow} (left). The function
$\tilde\varphi_\infty$ is discontinuous (cf. \cite[
(3.16)]{hmv1}),
\[
\tilde\varphi_\infty(\eta)=\cases{ 
 \frac{2C^d}{\eta^d} &
$\mbox{for}\quad \eta>C$\cr 0 &$\mbox{for}\quad \eta<C,$ }
\]
where $2C^d$ is the mass accumulated in the origin, which can be
chosen arbitrarily. In~\cite{hmv1} this type of blow-up was studied
using matched asymptotic expansions. There it was suggested that this
behaviour is stable and moreover it was expected that
there exist initial data such that (\ref{phity}) holds
uniformly in $\eta$ on compact subsets away from the shock.
A result of this type was proved in \cite[Theorem 3]{ehv} for a
related equation.
\section{Preliminaries}\label{prel}
\subsection{Estimates}
In this section we develop some estimates for  problem
(\ref{bequation})--(\ref{bequation2}),
which in turn will imply bounds for the  self-similar
problem (\ref{equaB1})--(\ref{equaB3}).
\begin{lemma}\label{leminic}
 If $b_0$ satisfies (\ref{inic}) then
\begin{eqnarray}
\label{est.br}
\frac r{d}b_r+b\geq 0\quad \mbox{in}\quad D_T.
\end{eqnarray}
\end{lemma}
{\bf Proof.} The solution $n$ of problem
(\ref{P1.s})--(\ref{mass1int.s}) satisfies the relation
\begin{eqnarray}\label{nintermb}
n=\frac 1{\chi_d}[\frac r{d}b_r+b]\quad \mbox{in}\quad D_\Tim.
\end{eqnarray}
Since $n_0\geq 0$ in $D,$  an application of the maximum principle
to problem (\ref{P1.s})--(\ref{mass1int.s}) shows that  $n\geq 0$
in $D_\Tim.$ Using this and (\ref{nintermb}) the result
follows.\eop

To prove the following results, we proceed as in \cite{fm} where
similar estimates were found for the semilinear parabolic equation
(\ref{semi-linear}).
\begin{lemma}
\label{lemma:br<0}
If $b_0$ satisfies (\ref{inicon0})
then
\begin{eqnarray}
\label{brneg}
b_r(r,t)<0 \quad\mbox{in}\quad D_T.
\end{eqnarray}
\end{lemma}
{\bf Proof.} Set $w(r,t):=r^{d+1}b_r(r,t).$
Differentiating (\ref{bequation}), we find
\begin{equation}
w_t-\chi_d\Theta\left(w_{rr}-\frac{d+1}{r}w_r\right)-\frac
1{d}rbw_r=\left(b+\frac 1{d}rb_r\right)w.
\end{equation}
Assume for the moment a stronger assumption on the initial data
\begin{eqnarray}\label{inicon}
(b_0)_r(r)<0\quad\mbox{for}\quad r\in (0,1)\quad \mbox{and}\quad
(b_0)_{rr}(0)<0.
\end{eqnarray}
This gives $w(0,r)=r^{d+1}b_r(0,r)< 0.$ Under (\ref{inicon}) the
function $b\equiv 1$ is a sub-solution for
(\ref{bequation})--(\ref{inic}), but not a solution; by Hopf's
Lemma, $w(1,t)=b_r(t,1)<0$ for all $t>0$, so that $w<0$ on $D_T$,
hence $b_r<0$ on $D_T$. To finish the proof, we note that by the
strong maximum principle, if  $b_0$ satisfies (\ref{inicon0}),
then for each $t_1\in (0,T)$ condition (\ref{inicon}) holds for
the function $b(r,t_1).$
This proves the result. \eop
\begin{lemma}
If $b_0$ satisfies  (\ref{inicon0}) and assuming that blow up
occurs at time $T>0$ , then
\begin{eqnarray}\label{bot-1}
b(0,t)\geq (T-t)^{-1}\quad \mbox{for}\quad t\in[0,T),
\end{eqnarray}
\end{lemma}
{\bf Proof.} Since the maximum of $b$ in $D$ is attained at $r=0$
(by Lemma~\ref{lemma:br<0}), we have $b_{rr}(0,t)\leq 0.$  It follows from
(\ref{bequation}) that $b_t(0,t)\leq b^2(0,t)$. Integrating this
inequality on $(0,T)$ gives  the result. \eop
\begin{lemma}\label{lemmabt}
If $b_0$ satisfies (\ref{btpositive})
then $b_t\geq 0$ for all $t\in (0,T)$.
\end{lemma}
{\bf Proof.} Condition (\ref{btpositive}) implies
that $b_0$ is a subsolution for~(\ref{bequation}-\ref{bequation2}); therefore
$b(r,\epsilon)\geq b(r,0)$ for small $\epsilon\geq0$. By
the comparison principle we find
$b(r,t+\epsilon)\geq b(r,t)$ for $t\in (0,T-\epsilon)$. It follows
that $b_t\geq0$ on $D_T$.\eop

The next lemma gives a  bound on $|b_r|$ in $D_T.$
\begin{lemma} Let $b_0$ satisfy (\ref{inicon0}) and (\ref{btpositive}).
Then
\begin{eqnarray}\label{bound3}
\chi_d\Theta b_r^2(r,t)\leq \frac 2{3}b(0,t)^3\quad
\mbox{for}\quad (r,t)\in D_T.
\end{eqnarray}
\end{lemma}
{\bf Proof.} Since  $b_t\geq 0$ and $b_r\leq 0$ in $D_T,$ we
multiply equation (\ref{bequation}) by $b_r$ and obtain
\begin{eqnarray*}
0&\geq &\chi_d \Theta\int\limits_0^rb_r b_{rr}\,ds+ \frac
1{3}b^3(r,t)-\frac 1{3}b^3(0,t) \\
&=& \frac12 \chi_d\Theta [b_r^2(r,t)-b_r^2(0,t)]+\frac 1{3}b^3(r,t)-\frac
1{3}b^3(0,t).
\end{eqnarray*}
Since $b_r^2(0,t)=0$ we obtain the desired inequality.\eop

To conclude this section we translate the properties of solutions
derived above into estimates for problem
(\ref{equaB1})--(\ref{equaB3}).
From hypothesis (\ref{corobB}) and noting that $b\geq 1$ and
$b_r\leq 0$ in $D_T,$
 we have the a priori bound
\begin{eqnarray}
\label{boundB} 0\leq B(\eta,\tau)\leq M\quad \mbox{for}\quad
(\eta,\tau)\in \Pi.
\end{eqnarray}
Combining this with (\ref{bound3})  and (\ref{brneg}),  we obtain
\begin{eqnarray}\label{boundBeta}
0\leq -B_\eta(\eta,\tau)\leq \bar M\quad\mbox{for}\quad
(\eta,\tau)\in \Pi,
\end{eqnarray}
where $\bar M$ depends on $M$. Finally from (\ref{bot-1}), we get
\begin{eqnarray}
\label{boundB0t} 1\leq B(0,\tau)
\quad \mbox{for}\quad \tau\in (0,\zmax).
\end{eqnarray}
\subsection{The steady state equation
(\ref{staW.n})}\label{sectionstaedystate}
 We begin by recalling
problem (\ref{staW.n}):
\begin{eqnarray}\label{staW2.a}
\varphi_{\eta\eta}+\frac {d+1}{\eta}\varphi_\eta+\frac 1{d}\eta
\varphi \varphi_\eta-
\frac 1{2}\eta \varphi_\eta+\varphi^2-\varphi=0 \quad\mbox{for}\quad \eta>0, \\
\varphi(0)\geq 1,\quad \varphi_\eta(0)=0. \label{staW2.b}
\end{eqnarray}
Condition (\ref{staW2.b}) is required, since $B(0,\tau)\geq 1$ for
all $\tau\geq 0.$
Equation  (\ref{staW2.a}) has three special solutions:
\[
\varphi_S(\eta)=\frac {2d}{\eta^2},\quad
\varphi^*(\eta)= 1,\quad\mbox{and}\quad \varphi_{*}(\eta)= 0\quad
\mbox{for}\quad \eta>0.
\]
Note that $\varphi_S$ satisfies
\begin{equation}\label{B_Sp}
\fl\qquad\varphi_S+\frac 1{2}\eta(\varphi_S)_\eta=0\quad
\mbox{and}\quad 0=(\varphi_S)_{\eta\eta}+\frac {d+1}{\eta}
(\varphi_S)_\eta+\frac 1{d}\eta
\varphi_S(\varphi_S)_\eta+(\varphi_S)^2.
\end{equation}
For bounded non-constant solutions we have the following
theorem~\cite{bcksv,hmv2}.
\begin{teo}\label{stasol}
Let $2<d<10.$ There exists a countable set of solutions
$\{\varphi_k\}_{k\in \N}$ of (\ref{staW2.a})--(\ref{staW2.b}) such
that $\varphi_k(0)>1$ and $\varphi_k(0)\to \infty$ as $k\to
\infty,$ Moreover $\varphi_k$ intersects the singular
solution $\varphi_S$ $k$ times   and has the asymptotic behaviour
$\varphi_k(\eta)\eta^2=Const(k)>0.$
\end{teo}
The proof is based on the equation for
$G(\eta):=\eta^2\varphi(\eta),$
\begin{eqnarray}\label{Gequation}
&G_{\eta\eta}+\left(\frac{(d-3)}{\eta}+\frac 1{d}\frac
G{\eta}-\frac 12\eta\right)
G_\eta+\frac{2(d-2)G}{\eta^2}\left(\frac G{2d}-1\right)=0,\\
\label{Gequation1} &\lim\limits_{\eta\downarrow 0} \frac
{G(\eta)}{\eta^2}<\infty,\quad \lim\limits_{\eta\to \infty}\eta
G_\eta(\eta)=0.
\end{eqnarray}
Note that $\varphi_S$ corresponds to $G(\eta)\equiv 2d.$

It was formally argued in~\cite{bcksv} that for each integer
$k\geq 2$ and $2<d<10$ the set
$$
{\cal S}_k=\{  \varphi \colon\mbox{$\varphi$ solution of
(\ref{staW2.a})--(\ref{staW2.b}) with $k$ intersections with
$\varphi_S$}\}
$$
is a singleton and that for $d>2$ the set ${\cal S}_1$
contains only two elements. More precisely,  ${\cal S}_1$ consists
of the functions $\varphi^*$ and $\varphi_1$ given in (\ref{phis1}).
If we relax condition (\ref{staW2.b}) to $\varphi(0)>0,$ we conjecture
that there is at least one other solution in ${\cal S}_1.$ For $d=3$ this was 
shown numerically by Brenner et al., who found a
solution $\varphi_1^*$ of (\ref{staW2.a}) such that
$\varphi_1^*(0)<1$ and $(\varphi_1^*)_\eta(0)=0,$ which intersects
$\varphi_S$ once~\cite[Figure 14]{bcksv}.
\section{Convergence}\label{convtheo}
In this section we prove the following convergence theorem.
\begin{teo}\label{converg}
Let conditions (\ref{inicon0}) and (\ref{btpositive}) hold. Let
$B(\eta,\tau)$ be a uniformly bounded global solution of
(\ref{equaB1})--(\ref{equaB3}). Then for every sequence $\tau_n\to
\infty$ there exists a subsequence $\tau_n'$ such that
$B(\eta,\tau_n')$ converges to a solution
$\varphi$ of (\ref{staW2.a})--(\ref{staW2.b}). The convergence is uniform on
every compact subset of $[0,\infty).$
\end{teo}
{\bf Proof.} Define $B^\sigma(\eta,\tau):=B(\eta,\sigma+\tau).$ We
will first show that for any unbounded
sequence $\{n_j\}$ there exists a subsequence (renamed
$\{n_j\}$) such that $B^{n_j}$ converges to a solution $\varphi$ of
(\ref{staW2.a})--(\ref{staW2.b}) uniformly in compact subsets of
$\R^+\times\R.$ Without loss of generality we assume that the sequence
$\{n_j\}$ is increasing.

Let $N\in \N$. We take $i$ large enough such that the rectangle
${\cal Q}_{2N}=\{(\eta,\tau)\in \R^2\colon 0\leq \eta\leq
2N,\:|\tau|\leq 2N\}$ lies in the domain of $B^{n_i}$. The function
$\tilde B(\xi,\tau)=B^{n_i}(|\xi|,\tau)$ is a solution of
$$
\tilde B_\tau =\Delta_{d+2} \tilde B- \frac 1{2}\xi\cdot\nabla
\tilde B+\frac 1{d}(\xi\cdot\nabla \tilde B)\tilde B+ \tilde
B^2-\tilde B
$$
on the cylinder given by
$$\Gamma_{2N}=\{(\xi,\tau)\colon \R^{d+2}\times \R\colon |\xi|\leq
2N,\:|\tau|\leq 2N\},
$$
and $|\tilde B(\xi,\tau)|$ is uniformly bounded in $\Gamma_{2N}$
by (\ref{boundB}).

By Schauder's interior estimates all partial derivatives of
$\tilde B$ can be uniformly bounded on the subcylinder
$\Gamma_N\subset \Gamma_{2N}.$ Consequently $B^{n_i},$
$B^{n_i}_\tau,$ $B^{n_i}_\eta,$ and $B^{n_i}_{\eta\eta}$ are
uniformly Lipschitz on ${\cal Q}_N\subset {\cal Q}_{2N}.$ By
Arzela-Ascoli, there is a subsequence $\{n_j\}_1^\infty$ and a
function $\bar B$ such that $B^{n_i},$ $B^{n_i}_\tau,$
$B^{n_i}_\eta,$ and $B^{n_i}_{\eta\eta}$ converge to $\bar B,$
$\bar B_\tau,$ $\bar B_\eta,$ and $\bar B_{\eta\eta}$, uniformly
on ${\cal Q}_N.$

Repeating the construction for all $N$ and taking a diagonal
subsequence, we can conclude that
\begin{eqnarray}\label{c.to.bar}
B^{n_j}\to \bar B, \quad B^{n_j}_\tau\to \bar B_\tau,\quad
B^{n_j}_\eta\to \bar B_\eta,\quad\mbox{and}\quad
B^{n_j}_{\eta\eta}\to \bar B_{\eta\eta},
\end{eqnarray}
uniformly in every compact subset in $\R^+\times \R.$
 Clearly $\bar B$ satisfies (\ref{equaB1}) and estimates (\ref{boundB})
 and (\ref{boundBeta}). Finally,
 it remains to prove that $\bar B$ is independent of $\tau.$
This implies that $\bar B$ is a  solution of (\ref{staW.n}), since
$B(0,\tau)\geq 1$ for all $\tau>0$, and the result follows.

{\bf Claim.} The function $\bar B$ is independent of $\tau.$

To prove this we construct a \emph{non-explicit} Lyapunov
functional in the spirit of Galaktionov \cite{ga1} and Zelenyak
\cite{zel}.

{\bf 1. Non-explicit Lyapunov functional.} We seek a Lyapunov
function of the form
$$
E(\tau)=\int\limits_0^\zmax\Phi(\eta,
B(\eta,\tau),B_\eta(\eta,\tau))\, d\eta,
$$
where $\zmax=(\chi_d\Theta T)^{-1/2}e^{\tau/2}$ and
$\Phi=\Phi(\eta,v,w)$ is a function to be determined.
In Appendix A we show that such a Lyapunov function exists;
more precisely, we show that a function $\rho=\rho(\eta,v,w)$ exists
such that
\begin{eqnarray}\nonumber
\frac d{d\tau}E(\tau)=-\int\limits_0^\zmax
\rho\big(\eta,B(\eta,\tau),B_\eta(\eta,\tau)\big)
(B_\tau)^2(\eta,\tau)\,d\eta \\
\qquad\qquad
+\left.\Phi_w B_\tau\right|_0^\zmax+\frac
1{2}\zmax\Phi\big(\zmax,B(\zmax,\tau),B_\eta(\zmax,\tau)\big).\label{lyapbound.n}
\end{eqnarray}
To identify the relevant domain of the functions $\Phi$ and $\rho$,
we note that by estimates
(\ref{boundB}) and (\ref{boundBeta}) the solution $B$ satisfies
$(\eta,B(\eta,\tau),B_\eta(\eta,\tau))\in \tilde \Rset,$ with
\begin{eqnarray}\label{Rtilde}
\tilde \Rset=\Rset\cap\{0\leq v\leq M,\: 0\leq -w\leq \bar M\},
\end{eqnarray}
where $\Rset=\{\eta>0,v\geq 0,w \leq 0\}\cup\{\eta=0,v\geq
0,w=0\}$.

The functions $\rho$ and $\Phi$ are continuous in
$\Rset\setminus\{\eta=\bar\eta,\: v>1\}$ with $\bar \eta>0$
defined later and they satisfy
\begin{eqnarray}\label{boundonrho}
\fl \qquad\frac
1{C_0}\eta^{d+1}e^{-C_0\eta^2}\leq\rho(\eta,v,w)\leq \eta^{d+1}
e^{-(d-2)\eta^2/4d}\quad\mbox{for}\quad (\eta,v,w)\in \tilde
\Rset,
\end{eqnarray}
with $C_0=C_0(M)>0$ (Lemma \ref{prho}), and
\begin{eqnarray}
\label{Phiest}
|\Phi(\eta,v,w)| \leq C_1\eta^{d+1} e^{-(d-2)\eta^2/4d}
\quad\mbox{for}\quad (\eta,v,w)\in \tilde \Rset
\end{eqnarray}
for some positive constants $C_1(M)>0$ (Lemma \ref{pPhi}).

{\bf 2. Proof of the claim.} An integration over the interval
$(a,b)$ of (\ref{lyapbound.n}) gives
\begin{eqnarray}\label{Epsi}
\fl\qquad\int\limits_a^b\int\limits_0^\zmax\rho(\eta,B(\eta,\tau),B_\eta(\eta,\tau))
B_\tau^2(\eta,\tau)\,d\eta d\tau=E(a)-E(b)+\psi(a,b)
\end{eqnarray}
where
\begin{eqnarray}\nonumber
\fl\qquad\psi(a,b)&:=&\int\limits_a^b \frac
1{2}\zmax\Phi(\zmax,B(\zmax,\tau),B_\eta(\zmax,\tau))\, d\tau + \\
\label{psi1}
\fl\qquad &&+\int\limits_a^b
B_\tau(\zmax,\tau)\left[\int\limits_0^{B_\eta(\zmax,\tau)}\rho(\zmax,B(\zmax,\tau),s)\,ds\right]d\tau.
\end{eqnarray}
Since $B_\tau(\zmax,\tau)=-B(\zmax,\tau)-\frac 1{2}\zmax
B_\eta(\zmax,\tau),$
\[
B_\tau(\zmax,\tau)=-Te^{-\tau}-\frac 1{2}b_r(1,T(1-e^{\tau})).
\]
Applying~(\ref{est.br}) at $r=1$ gives
$|b_r(1,T(1-e^{\tau}))|\leq d$ and consequently $B_\tau$ is
uniformly bounded as $\tau\to\infty$. Employing this bound on
$B_\tau$ and the estimates (\ref{boundonrho}) and
(\ref{Phiest}) we find
\begin{eqnarray}\label{limphi}
\lim\limits_{a\to\infty}\{\sup\limits_{b>a}\psi(a,b)\}=0.
\end{eqnarray}
By (\ref{c.to.bar}), we have that there exists a sequence $n_j\to
\infty$ such that $B^{n_j}(\eta,\tau)$ converges to $\bar B$
uniformly in compact subsets of $(\R^+)^2.$ For any fixed $N$ we
will prove for a subsequence satisfying
$\lim\limits_{j\to\infty}(n_{j+1}-n_j)=\infty$ that
\begin{eqnarray}
\label{Bnj.to.0}
\lim\limits_{n_j\to\infty}\int\limits_{{\cal
Q}_N}\rho(\eta,B^{n_j}(\eta,\tau),B_\eta^{n_j}(\eta,\tau))
(B_\tau^{n_j})^2(\eta,\tau)\,d\eta d\tau=0,
\end{eqnarray}
where we recall that
${\cal Q}_{N}=\{(\eta,\tau)\colon \R^2\colon 0\leq \eta\leq
N,\:|\tau|\leq N\}.$
Since $\rho$ is bounded from below on bounded subsets of
$\tilde\Rset$, it then follows that
\[
\int\limits_{{\cal Q}_N} \bar B_\tau^2 \, d\eta d\tau =
\lim\limits_{n_j\to\infty}
  \int\limits_{{\cal Q}_N}(B_\tau^{n_j})^2(\eta,\tau)\,d\eta d\tau=0,
\]
proving the claim.
For all $j$ sufficiently large,
$$
N\leq (\chi_d\Theta T)^{-1/2}e^{\frac
1{2}(n_j-N)}\quad\mbox{and}\quad n_{j+1}-n_j\geq 2N.
$$
Consequently using (\ref{Epsi}), we find
\begin{eqnarray*}
\fl\int\limits_{-N}^N\int\limits_0^N\rho(\eta,B^{n_j}(\eta,\tau),B_\eta^{n_j}(\eta,\tau))
(B_\tau^{n_j})^2(\eta,\tau)\,d\eta d\tau \\
\fl\qquad\leq
\int\limits_{-N}^{-N+n_{j+1}-n_j}\int_0^{(\chi\Theta^*T)^{-1/2}
e^{\frac1{2}(n_j-N)}}\rho(\eta,B^{n_j}(\eta,\tau),B^{n_j}_\eta(\eta,\tau))
(B_\tau^{n_j})^2(\eta,\tau)\,d\eta d\tau \\
\fl\qquad\leq
\int\limits_{n_j-N}^{n_{j+1}-N}\int_0^{(\chi\Theta^*T)^{-1/2}
e^{\frac1{2}(n_j-N)}}\rho(\eta,B(\eta,\tau),B_\eta(\eta,\tau))
(B_\tau)^2(\eta,\tau)\,d\eta d\tau \\
\fl\qquad\leq E(n_j-N)-E(n_{j+1}-N)+\psi(n_j-N,n_{j+1}-N).
\end{eqnarray*}
Hence applying (\ref{limphi}), we discover
\[
\fl\int\limits_{{\cal
Q}_N}\rho(\eta,B^{n_j}(\eta,\tau),B_\eta^{n_j}(\eta,\tau))
(B_\tau^{n_j})^2(\eta,\tau)\,d\eta d\tau\leq
\limsup\limits_{j\to\infty}[E(n_j-N)-E(n_{j+1}-N)].
\]
Next we divide the expression $E(n_j-N)-E(n_{j+1}-N)$ into three
integrals, choosing $K$ arbitrarily large:
\begin{eqnarray}
\nonumber
\fl E(n_j-N)-E(n_{j+1}-N)=\\
\fl\qquad\label{Eexp1}=\int\limits_0^K[\Phi(\eta,B^{n_j}(\eta,-N),B^{n_j}_\eta(\eta,-N))
      -\Phi(\eta,B^{n_j}(\eta,-N),B^{n_j}_\eta(\eta,-N)]\,d\eta\\
\fl\qquad\label{Eexp2}
  +\int\limits_K^{T^{-1/2}e^{\frac{n_j-N}{2}}}
\Phi(\eta,B^{n_{j+1}}(\eta,-N),B^{n_{j+1}}_\eta(\eta,-N))\,d\eta\\
\fl\qquad\label{Eexp3}
  +\int\limits_K^{T^{-1/2}e^{\frac{n_{j+1}-N}{2}}}
     \Phi(\eta,B^{n_j}(\eta,-N),B^{n_j}_\eta(\eta,-N))\,d\eta.
\end{eqnarray}
Integral (\ref{Eexp1}) tends to zero as $j\to \infty.$ In fact by
the continuity of $\Phi$ in the second and third argument we
obtain  pointwise convergence and by the bounds (\ref{Phiest}) on
$\Phi,$ we apply the Dominated Convergence Theorem to conclude.
Expressions (\ref{Eexp2}) and (\ref{Eexp3}) can be made
arbitrarily small since they can be bounded by
$$
C\int\limits_K^\infty\eta^{d+1}e^{-(d-2)\eta^2/4d}d\eta,
$$
where $C$ is a positive constant, and $K$ can be chosen arbitrary large.
Thus we have proved~(\ref{Bnj.to.0}), concluding the proof of
the Theorem.
\eop

\section{Comparison results}\label{intersectioncomparison}
\subsection{Comparison with the singular solution $\varphi_S$}

This section closely follows~\cite{be}.
From section \ref{sectionstaedystate}, we recall that solutions
$\varphi$ of (\ref{staW2.a})--(\ref{staW2.b}) are classified by
their intersections with $\varphi_S.$ In this section we study the
intersections of solutions $B$ of (\ref{equaB1})--(\ref{equaB3})
with $\varphi_S.$ Our results are closely related to the ones
found in \cite{be}, where equation (\ref{semi-linear}) was
studied.

We first see that for $\Theta< 1/(2d\chi_d)$ a solution $B$ of
(\ref{equaB1})--(\ref{equaB3}) intersects the singular solution
$\varphi_S$ at least once in $\Pi(0)$ since
$$
\varphi_S(0)=\infty>B(0,0),\:\mbox{and}\:\varphi_S\left((\chi_d\Theta
T)^{-1/2}\right)<B\left((\chi_d\Theta T)^{-1/2},0\right)=T.
$$
On the other hand, for $\Theta\geq 1/(2d\chi_d)$ it can also be
shown that $B$ intersects $\varphi_S$ at least once in $\Pi(0).$
Assuming the contrary, suppose that $B(\cdot,0)<\varphi_S(\cdot)$
in $\Pi(0).$  By the maximum principle, we obtain $B<\varphi_S$ in
$\Pi.$ Therefore in the limit $\tau\to \infty,$  thanks to the
Theorem \ref{converg} and since $B(0,\tau)\geq 1$ for all
$\tau>0$, we find a solution $\varphi$ of (\ref{staW.n}) such that
$\varphi <\varphi_S.$ However we can show  that every bounded non
zero solution $\varphi$ of (\ref{staW.n}) has to cross
$\varphi_S.$
This is equivalent to proving that there exists no solution $G$ of
(\ref{Gequation})--(\ref{Gequation1}) such that $G(\eta)<2d$ for
$\eta\geq 0.$ To check this, we assume that such a solution exists;
we examine two cases. Suppose that for some
$\eta^*$, we have $G_\eta(\eta^*)=0$ and $G(\eta^*)<2d.$ By (\ref{Gequation}), $G$
has a strict minimum at $\eta^*$, which contradicts the  boundary
condition~(\ref{Gequation1}). On the
other hand if, $G(\eta)$ is increasing for all $\eta>0$, then for large
$\eta,$ equation (\ref{Gequation}) implies that $G_{\eta\eta}>0,$
which also contradicts (\ref{Gequation1}).

We conclude that there exists $\eta_1\in \Pi(0)$ such that
$B(\eta_1,0)=\varphi_S(\eta_1)$ and $B(\eta,0)<\varphi_S(\eta)$
for $\eta<\eta_1.$
\begin{lemma}\label{confun}
Under the assumptions (\ref{inicon0}) and  (\ref{btpositive}),
there exists a con\-ti\-nuously di\-ffe\-rentiable function
$\eta_1(\tau)$ with domain $[0,\infty)$ such that
$\eta_1(0)=\eta_1$ and
$B(\eta_1(\tau),\tau)=\varphi_S(\eta_1(\tau))$ for all $\tau\geq
0.$
\end{lemma}
{\bf Proof.} Define $H(\eta,\tau):=B(\eta,\tau)-\varphi_S(\eta).$
We first claim that $H,H_\eta,$ and $H_\tau$ do not vanish
simultaneously. Using Lemma \ref{lemmabt} and the strong maximum
principle we find
\begin{eqnarray}\label{btp.b}
b_t=(T-t)^{-2}\left(B_\tau+B+\frac 1{2}\eta B_\eta\right)>0\quad
\mbox{in}\quad D_T.
\end{eqnarray}
Suppose there exists a point in $\Pi$ where $H_\eta=H_\tau=H=0.$
Then $H_\tau=0$ implies $B_\tau=0,$ and  condition  $H_\eta=0$
combined with $H=0$ gives
$$
B+\frac 1{2}\eta B_\eta=0\quad \mbox{in}\quad \Pi,
$$
using  (\ref{B_Sp}). This implies that $b_t=0$ at some point of
$D_T,$ a contradiction with (\ref{btp.b}). Secondly, we claim that
$H_\eta\neq 0$ at any point $(\bar \eta,\bar \tau)\in \Pi$ where
$H(\bar \eta,\bar \tau)=0$ and moreover $H(\eta,\bar \tau)<0$ in a
left neighborhood of $\bar \eta.$
A proof of this can be done as in~\cite{be}. Moreover, from the
proof, we find $H_\eta(\bar \eta,\bar \tau)>0.$

Now we prove that $H_\eta(\eta_1,0)>0.$ This follows from the
equation satisfied by $H(\eta,0).$
To the left of $\eta_1,$ we find
\begin{eqnarray}
H_{\eta\eta}(\eta,0)&+\frac {d+1}{\eta}H_\eta(\eta,0)+\frac
1{2d}\eta
H_\eta(\eta,0)(B(\eta,0)+\varphi_S)\\
&+\frac 1{2d}\eta H(\eta,0)(B(\eta,0)+\varphi_S)_\eta\geq 0.
\end{eqnarray}
Since $(B(\eta,0)+\varphi_S)_\eta\leq 0$ and $H(\eta_1,0)=0,$ we
can apply Hopf's Lemma to obtain that $H_\eta(\eta_1,0)>0.$%
Finally, to conclude the proof of the lemma, we use
the implicit function theorem as in~\cite{be}.\eop

Define the set $ \Pi_1=\{ (\eta,\tau) \mid 0<\eta<\eta_1(\tau)\:\}
$ and the function
$$
\eta_2(\tau)=e^{\tau/2}\cdot
  \sup\Bigl\{\eta\in (\eta_1,(\chi_d\Theta T)^{-1/2}]\colon
H(s,0)\geq 0\:\mbox{for}\:s\in [\eta_1,\eta]\Bigr\}.
$$
Since $H(\eta_1,0)=0$ and $H_\eta(\eta_1,0)>0,$ the above supremum
is finite. Define the set
\[
\Pi_2=\{ (\eta,\tau) \mid \eta_1(\tau)<\eta<\eta_2(\tau)\}. \] Let
$F(\tau)=H(\eta_2(\tau),\tau).$ By definition of $\eta_2,$
$F(0)\geq 0.$ Also, \[ \frac
d{d\tau}F(\tau)=H_\tau(\eta_2(\tau),\tau)+\frac
1{2}\eta_2(\tau)H_\eta(\eta_2(\tau),\tau).
\]
Using (\ref{btp.b}), we have $ d[e^\tau F(\tau)]/{d\tau}\geq 0.$
An integration yields $F(\tau)\geq 0$ for $\tau\geq 0.$

As  was done in \cite{be}, applying  the maximum principle, using
Lemma \ref{confun}, and noting that $H(\eta_2(\tau),\tau)\geq 0$
for $\tau\geq 0$, we can prove the following lemma and its
corollary.
\begin{lemma} The function $H(\eta,\tau)=B(\eta,\tau)-\varphi_S(\eta)$ satisfies $H<0$ in $\Pi_1$ and $H>0$ in $\Pi_2.$
\end{lemma}
\begin{coro}\label{intsing}
Assume the conditions in Lemma \ref{confun}. For each $N>0$ there
is $\tau_N>0$ such that for $\tau>\tau_N,$ $B(\eta,\tau)$
intersect $\varphi_S(\eta)$ at most once in $\eta\in (0,N).$
\end{coro}
\subsection{Intersection comparison}
 In this section we derive  comparison results,  which
will be used to prove that $\omega,$ the  limit set
 (\ref{omegalimitset}),
is a singleton.

We start by considering  the following linear equation with
inhomogeneous boundary conditions:
\begin{equation}\label{vequation}
\fl\quad\cases{ 
v_t=v_{rr}+\frac{d+1}{r}v_r+a(r,t)v,  & for $0<r<1,\:T_1<t<T_2,$
\cr
v_r(0,t)=0 \:& for $T_1<t<T_2,$ \cr v(1,t)=h(t)\:& for $T_1<t<T_2;$ }
\end{equation}
where $T_1,T_2$ are positive constants and
\begin{equation}\label{hypvequation}
a\in L^\infty([0,1]\times(T_1,T_2)),\quad h\in C^1((T_1,T_2)),
\end{equation}
are given functions. Moreover we assume
\begin{equation}\label{hypvequation1}
h(t)>0\quad \mbox{for}\quad T_1<t<T_2.
\end{equation}
The zero number functional of (\ref{vequation}) is
defined by
\begin{equation}
z[v(\cdot,t)]=\#\{r\in[0,1]\colon v(r,t)=0\},
\end{equation}
and the following lemma  provides some properties of this zero number
functional.
\begin{lemma}[\cite{m}]\label{Matoslemma}
Let $v=v(r,t)$ be a nontrivial classical solution of
(\ref{vequation}) and assume that (\ref{hypvequation}) and
(\ref{hypvequation1}) hold. Then the following properties hold
true:
\begin{enumerate}
\item $z[v(\cdot,t)]<\infty$ for any $T_1<t<T_2;$ \item
$z[v(\cdot,t)]$ is nonincreasing in time; \item if
$v(r_0,t_0)=v_r(r_0,t_0)=0$  for some $r_0\in [0,1]$ and
$t_0>T_1,$ then $z[v(\cdot,t)]$ drops strictly at $t=t_0,$ that
is, $z[v(\cdot,t_1)]>z[v(\cdot,t_2)]$ for any
$T_1<t_1<t_0<t_2<T_2.$
\end{enumerate}
\end{lemma}
From this lemma we deduce a property of intersection between a
solution $\varphi$  of (\ref{staW.n})  and a solution $B$  of
(\ref{equaB1})--(\ref{equaB3}).

\begin{lemma}\label{ourlemma}
Let $B$ be a bounded solution of (\ref{equaB1})--(\ref{equaB3})
and let $\varphi$ be a solution of (\ref{staW2.a}).
Denote $Z(\tau)=\#\{r\in[0,\zmax]\colon
B(\eta,\tau)=\varphi(\eta)\}.$ Then the following properties hold
true:
\begin{enumerate}
\item $Z(\tau)<\infty$ for any $\tau>\tau^*;$
\item $Z(\tau)$ is
nonincreasing in time;
\item if $B(\eta_0,\tau_0)=\varphi(\eta_0)$ and
$B_\eta(\eta_0,\tau_0)=\varphi_\eta(\eta_0)$ for $\tau_0>\tau_1,$
and $\eta_0\leq \zmaxa[\tau]$ then  $Z(\tau_1)>Z(\tau_2)$ for any
$\tau_1<\tau_0<\tau_2.$
\end{enumerate}
\end{lemma}

{\bf Proof.} Writing $\bar V=U-b$, where
$U(r,t)=(T-t)^{-1}\varphi(r/(\chi_d\Theta(T-t))^{1/2})$, we have
\begin{eqnarray}
\label{difequation}
\fl\cases{ 
\bar V_t=\bar V_{rr}+\left(\frac{d+1}{r}+\frac r{d}U\right)\bar
V_r+\left(\frac r{d}b_r+b+U\right)\bar V & $\mbox{for}\: 0<r<1,\:0<t<T,$ \cr
\bar V_r(0,t)=0,\quad \bar V(1,t)=U(1,t)-b(1,t) & $\mbox{for}\: 0<t<T.$
}
\end{eqnarray}
Let $T_1<T_2<T.$ For the variable  $\dis
V(r,t)=\exp\left(\frac 1{2d}\int\limits_0^r yU(y,t)\,dy\right)\bar
V(r,t),$ we find
\begin{eqnarray*} 
\fl\cases{ 
 V_t=V_{rr}+\frac{d+1}{r} V_r+A(r,t) V & $\mbox{for}\quad
0<r<1,\:T_1<t<T_2,$\cr 
 V_r(0,t)=0, & $\mbox{for}\quad T_1<t<T_2,$ \cr
 V(1,t)=(U(1,t)-1)\exp\left(\frac
1{2d}\int\limits_0^1
yU(y,t)\,dy\right) & $\mbox{for}\quad T_1<t<T_2,$
}
\end{eqnarray*}
where
\[
\fl
A(r,t) = \frac rd b_r +b + U + \frac1{2d} \int_0^r yU_t(y,t)\, dy
  - \frac1{4d^2} r^2U^2 - \frac1{2d}(U+rU_r) - \frac{d+1}{2d}U.
\]
Note that  $A\in L^\infty([0,1]\times(T_1,T_2))$ since
$b,b_r,U,U_t,U_r\in L^\infty([0,1]\times(T_1,T_2)).$ If we show that
$V(1,t)$ does not change sign for $t>t_0,$ then setting $T_1=t_0$ and using
Lemma \ref{Matoslemma}, we have proved the lemma.

We claim that there exists $\bar t_0$ such that
$U_t(1,t)$ does not change sign for $t>\bar t_0.$  By definition of $V,$ this implies that there exists $t_0\geq \bar t_0$ such that $V(1,t)$ does not change
sign for $t>t_0.$

Since $U_t(r,t)=(T-t)^{-2}(\eta^2\varphi)_\eta/(2\eta),$ if $r=1$ and $t>t^*,$ then
\begin{equation}\label{ut1t}
\fl \qquad U_t(1,t)=(T-t)^{-2}\frac 1{2\eta}(\eta^2\varphi)_\eta\quad \mbox{for}\quad t>t^*,\quad
\mbox{and}\quad
\eta>\eta^*(t^*),
\end{equation}
where $\eta^*(t^*):=(\chi_d\Theta(T-t^*))^{-1/2}.$
From  \cite[Lemma A.1]{bcksv0}, we know that for a given $a\in(0,4d),$ 
any solution $\varphi$ of  (\ref{staW2.a}) satisfying
\begin{equation}\label{asyb}
\eta^{2}\varphi(\eta)\to a\quad\mbox{as}\quad \eta\to\infty,
\end{equation}
is such that there exists $\bar\eta_0=\bar\eta_0(a)$ so that the sign of  $(\eta^2\varphi)_\eta$ 
does not change on $[\bar\eta_0,\infty).$
Using (\ref{ut1t}), this implies that there exists $\bar t_0=\bar t_0(\bar \eta_0)$ such that the claim holds.
 \eop

\section{Proofs of main results}
\label{proofresults}
We start by proving that the $\omega$-limit set of  problem
(\ref{equaB1})--(\ref{equaB3}) is a singleton.
\begin{theorem}
Assume the hypotheses of Theorem \ref{prin.s}. Then the  set
$\omega$ defined in (\ref{omegalimitset}) is a singleton.
\end{theorem}
{\bf Proof.} For this proof we extend a solution  $B$ of
(\ref{equaB1})--(\ref{equaB3}) to all $(\R^+)^2$ by setting
$B(\eta,\tau)=e^{-\tau}T$ for $(\eta,\tau)\in
(\R^+)^2\setminus\Pi.$ We also define the weight function
$\rho^*(\eta)=e^{-\eta^2/4}$ for $\eta>0$.

The hypothesis~(\ref{corobB}) implies that $B$ is uniformly bounded;
Theorem~\ref{converg} therefore states that $\omega$ is non-empty, and
that each $\varphi\in\omega$ is a solution of~(\ref{staW2.a})--(\ref{staW2.b}).

We claim that for each $\varphi\in \omega$ there exists
$\tau^*>0$ such that $B(0,\tau)-\varphi(0)$ never changes sign in
$[\tau^*,\infty).$ By contradiction, we assume that there exists a
sequence $\tau_k,$ such that $\tau_k\to \infty,$ and
$B(0,\tau_k)=\varphi(0).$ Since
$B_\eta(0,\tau_k)=\varphi_\eta(0)=0,$ by Lemma \ref{ourlemma} the
function $Z(\tau)$ has to decrease at least by one. However this
cannot happen an infinite number of times.
This proves the claim.

Suppose now that $\omega$ is not a singleton. Since the
$\omega$-limit set is connected, closed, and non empty,
it contains an infinite number of elements.
We select three different elements
$\varphi_1,\varphi_2,\varphi_3$ in the $\omega$-limit set. Since
these functions are different and each solves (\ref{staW.n}), we may  assume that
$\varphi_1(0)<\varphi_2(0)<\varphi_3(0).$ By the  claim above,
$B(0,\tau)-\varphi_2(0)$ never changes sign in $[\tau^*,\infty).$
This contradicts the fact that $\varphi_1$ and $\varphi_3$ are
elements of $\omega$; it follows that $\omega$ is a
singleton.\eop
\medskip

We now conclude the proof of Theorems \ref{prin.s} and  \ref{prin.s2}, and
Corollary \ref{prin.s3}.
\smallskip

{\bf Proof of Theorem \ref{prin.s}.} By
the previous theorem $\omega$ is a singleton, say $\{\bar B\}.$
From Corollary \ref{intsing}, we find that for every $N>0$ there
exists a $\tau_N>0$ such that the solution $B(\eta,\tau)$ intersects
$\varphi_S(\eta)$ at most once in $\eta\in[0,N]$ for each
$\tau>\tau_N.$ This implies that in the limit $\tau\to\infty,$
$\bar B$ intersects $\varphi_S$ at most once, concluding the
proof.\eop
\smallskip

{\bf Proof of Theorem \ref{prin.s2}.} Since  $b$ and
$U_1(r,t)=(T-t)^{-1}\varphi_1(r/(\chi_d\Theta(T-t))^{1/2})$ are
solutions of (\ref{bequation}) with the same blow up time, $\bar
V=b-U_1$ satisfies equation (\ref{difequation}). Using that
$U_1(r,t)=2d/((d-2)(T-t)+r^2/(2\chi_d\Theta))$, we find
$$
\bar V(1,t)=(1-U_1(1,t))>0\quad \mbox{if}\quad \Theta\leq
1/(4d\chi_d)\quad \mbox{for any}\quad t<T.
$$
The functions $U_1$ with $b$ necessarily intersect exactly once
for all $t$, since non-intersection implies that the solutions must have
different times of blow-up~\cite[p. 271]{sgkm}. It follows that
$b(0,0)< U_1(0),$ and one finds $(T-t)b(0,t)\leq 2d/(d-2)$.
An application of Theorem \ref{prin.s} proves
the theorem.\eop
\smallskip

{\bf Proof of Corollary \ref{prin.s3}.} If $b_0\equiv 1$ and
$\Theta< 1/(2(d+2)\chi_d)$, we know from \cite[Theorem 2]{bn1}
that the corresponding solution $b$ blows up. Now assuming
$\Theta\leq 1/(4d\chi_d)<1/(2(d+2)\chi_d),$ we can apply Theorem
\ref{prin.s2} to conclude.\eop

\appendix
\def\thesection{\Alph{section}}

\section{Appendix: The Lyapunov functional}
\label{sec:appendix}
In this appendix we construct the Lyapunov functional $E$
satisfying (\ref{lyapbound.n}) with  the suitable  properties of
$\rho$ and $\Phi$ to prove  Theorem \ref{converg}. We start with a
formal construction  of the functional.
This requires solving  a first-order equation for $\rho$
after which $\Phi$ can be  expressed  in terms of $\rho$.
Finally, we explain how to use smooth approximations of $\Phi$ to obtain a
rigorous derivation of (\ref{lyapbound.n}).
\subsection{Formal derivation of  a Lyapunov functional}
Assume that $\Phi$ and $\rho$ are regular. To find such
functions  satisfying (\ref{lyapbound.n}), we compute
\begin{eqnarray}\label{apen1}
\fl\frac d{d\tau}E(\tau)=\int\limits_0^\zmax\Phi_v B_\tau\,
d\eta+\int\limits_0^\zmax\Phi_w B_{\tau\eta}\, d\eta +\frac
\zmax{2}\Phi(\zmax,B(\zmax,\tau),B_\eta(\zmax,\tau)).
\end{eqnarray}
Wherever possible we omit the arguments of $\Phi$ and $\rho$ for clarity.
Integrating by parts the second integral in (\ref{apen1}) becomes
\begin{eqnarray}\nonumber
\int\limits_0^\zmax\Phi_w B_{\tau\eta}\,
d\eta=
\label{apen2}-\int\limits_0^\zmax[\Phi_{\eta
w}+\Phi_{vw}B_\eta+\Phi_{ww}B_{\eta\eta}] B_{\tau}\,
d\eta+\left.\Phi_w B_\tau\right|_0^\zmax.
\end{eqnarray}
Defining
\[
f(\eta,v,w)=\frac{d+1}{\eta}w-\frac{\eta}{2}w+\frac 1{d}\eta v
w+v^2-v,
\]
equation (\ref{equaB1}) takes the form
$B_\tau=B_{\eta\eta}+f(\eta,B,B_\eta)$,
by which equation~(\ref{apen1}) becomes
\begin{eqnarray*}\label{Phiformula}
\fl\qquad\frac d{d\tau}E(\tau)&=\int\limits_0^\zmax\left \{
[\Phi_v-\Phi_{\eta w}-\Phi_{vw}B_\eta+\Phi_{ww}f] B_{\tau}-
\Phi_{ww}(B_\tau)^2\right\}\,d\eta\\
&+\left.\Phi_w B_\tau\right|_0^\zmax+\frac
\zmax{2}\Phi(\zmax,B(\zmax,\tau),B_\eta(\zmax,\tau)).
\end{eqnarray*}
Now if functions $\rho=\rho(\eta,v,w)>0$ and $\Phi=\Phi(\eta,v,w)$
exist that satisfy the system of equations
\begin{eqnarray}\label{Phicond}
-\Phi_v+\Phi_{\eta w}+w\Phi_{vw}=\rho f\quad\mbox{and}\quad
\Phi_{ww}=\rho,
\end{eqnarray}
then  $E$ has the form of a  Lyapunov
functional with a contribution on the boundary, i.e.
\begin{eqnarray}\nonumber
\fl\qquad\frac d{d\tau}E(\tau)=&-\int\limits_0^\zmax
\rho(\eta,B,B_\eta) (B_\tau)^2\,d\eta  +\left.\Phi_w
B_\tau\right|_0^\zmax+
\\
 &+\frac
\zmax{2}\Phi(\zmax,B(\zmax,\tau),B_\eta(\zmax,\tau)).\label{lyapbound}
\end{eqnarray}
Therefore we may obtain this formula by solving system
(\ref{Phicond}), which we do by transforming it to  a first-order
equation for $\rho,$
\begin{eqnarray}\label{f.ord.rho}
w\rho_v+\rho_\eta-f\rho_w=f_w\rho.
\end{eqnarray}
If we supplement a given solution $\rho$ of this equation with the
function $\Phi$ given by
\begin{eqnarray}\label{pPhi2}
\Phi(\eta,v,w)=\int\limits_0^w (w-s)\rho(\eta,v,s)\,ds
-\int\limits_0^v\rho(\eta,\mu,0)f(\eta,\mu,0)\,d\mu,
\end{eqnarray}
then the pair $(\rho,\Phi)$ solves~(\ref{Phicond}). In order
to find the pair $(\rho,\Phi)$ we therefore only need to solve
equation~(\ref{f.ord.rho}).
\subsection{The first-order equation for $\rho$}
We solve equation (\ref{f.ord.rho}) by the method of
characteristics.
Characteristic curves of equation (\ref{f.ord.rho}) are curves
${\mathbf x}=(\eta,v,w)$ in $\R^3$, which we consider
parametrised by $\eta$, along which
\begin{eqnarray}\label{firstor}
\frac d{d\eta} v=w\qquad \mbox{and}\qquad
\frac d{d\eta} w=-f.
\end{eqnarray}
If a curve ${\mathbf x}(\eta)=(\eta,v^1(\eta),w^1(\eta))$ satisfies
these equations,
then equation~(\ref{f.ord.rho}) reduces
to
\begin{eqnarray}\label{firstor1}
\frac d{d\eta} \rho(\mathbf x(\eta)) = f_w(\mathbf x(\eta))  \rho(\mathbf x(\eta)).
\end{eqnarray}
In order to solve the system of ODE's (\ref{firstor}) and (\ref{firstor1}),
we select a vector
$(\eta_0,v_0,w_0)\in \R^+\times\R^2$ and define
$\phi(\xi)=\phi(\xi;\eta_0,v_0,w_0) $ to be the solution
of the initial value problem
\begin{eqnarray}\label{phieq}
\phi''+f(\xi,\phi,\phi')=0,\quad\mbox{with}\quad
\phi|_{\xi=\eta_0}=v_0\quad \mbox{and}\quad
\phi'|_{\xi=\eta_0}=w_0,
\end{eqnarray}
where $'=\frac \partial{\partial \xi}$. If the
curve $\mathbf x$ passes through $(\eta_0,v_0,w_0)$, i.e.
if $\mathbf x(\eta_0) = (\eta_0,v_0,w_0)$, then this curve can  be
identified with $\phi(\cdot;\eta_0,v_0,w_0)$, since
$\mathbf x(\eta) = (\eta,v^1(\eta),w^1(\eta))$ where
\begin{eqnarray}\label{forvw}
v^1(\eta)=\phi(\eta;\eta_0,v_0,w_0)\quad \mbox{and}\quad
w^1(\eta)=\phi'(\eta;\eta_0,v_0,w_0).
\end{eqnarray}
Since
$f_w=\frac {d+1}{\eta}-\frac{\eta}{2}+\frac 1{d}\eta v$,
we may integrate (\ref{firstor1}) to find
\begin{eqnarray}
\nonumber
\fl\qquad\rho(\eta,v,w)&=&\rho(\eta_0,v_0,w_0)\exp\left\{
\int\limits_{\eta_0}^\eta \left[\frac{d+1}{\xi}-\frac{\xi}{2}+\frac
1{d}\xi v^1(\xi)\right]\, d\xi\right\} \\
\fl\qquad &=& \rho(\eta_0,v_0,w_0) \frac{\eta^{d+1}}{\eta_0^{d+1}}
  e^{-\eta^2/4 + \eta_0^2/4}
  \exp \left\{\frac1d \int\limits_{\eta_0}^\eta \xi v^1(\xi)\, d\xi\right\}.
\label{def:rho}
\end{eqnarray}
To prove Theorem \ref{converg}, we need to define $\rho$
in the set $\tilde \Rset\subset
\Rset$ given by (\ref{Rtilde}),
\begin{eqnarray*}
\Rset=\{\eta>0,v\geq 0,w \leq 0\}\cup\{\eta=0,v\geq
0,w=0\} \\
\tilde \Rset=\Rset\cap\{0\leq v\leq M,\: 0\leq -w\leq \bar M\}.
\end{eqnarray*}
We do so in the following way: for each $(\eta,v,w)\in \Rset,$
we define $\rho(\eta,v,w)$ by
following the characteristic curve through $(\eta,v,w)$ to a
reference point $(\eta_0,v_0,w_0)$ for which
$\rho(\eta_0,v_0,w_0)$ is fixed by choice; the value of $\rho(\eta,v,w)$
is then given by~(\ref{def:rho}). To select an appropriate
set of reference points, we study some of the properties of
solutions $\phi$ of (\ref{phieq}), since they define the
characteristic curves.

It follows  from standard ODE theory that solutions of
(\ref{phieq}) are locally smooth and continuous under
changes of $(\eta_0,v_0,w_0)$.
In general, however, we cannot extend
these solutions to the whole of $\R^+$; in fact, for each
$(\eta,v,w)\in\Rset$, there  may exist $0\leq \xi_1 < \eta$
and/or $\xi_2>\eta$
such that
\[
\phi(\xi_1;\eta,v,w)=\infty\qquad
\mbox{and/or}\qquad\phi(\xi_2;\eta,v,w)=-\infty.
\]
Partly because of this difficulty, we choose to only use forward solutions of
(\ref{phieq}) to define the characteristic curves. The next result
details the behaviour of a forward solution $\phi$ of~(\ref{phieq}).
\begin{lemma}\label{claphi}
Let $(\eta,v,w)\in \Rset,$ and let
$\phi(\xi)=\phi(\xi;\eta,v,w)$ be the solution of (\ref{phieq}).
For $\xi\geq \eta,$ exactly one of the following three alternatives holds:
\begin{enumerate}
\item $\phi\equiv 1 $  or $\phi\equiv 0$;
\item there exists $\eta^*>\eta$
such that $\phi(\eta^*)=0$ and $\phi(\xi)<0$ for $\xi>\eta^*;$
\item $\phi(\xi)\to 0$ as $\xi \to \infty$ and  there exists a
constant $C>0$ such that $\phi(\xi)\xi^{2}\to C$ as $\xi\to
\infty$.
\end{enumerate}
\end{lemma}
{\bf Proof.} See \cite[p. 95]{g}. The proof is  based on results
from~\cite{hmv2}. \eop 
Since we need to define $\rho$ with
the appropriate estimates, we introduce a parameter~$\bar \eta$ in
the following lemma.

\begin{lemma}
\label{lemetabar2}
There exists $\bar\eta>0$ such that for every $\eta_1\geq \bar \eta$ any solution~$\phi$
of~(\ref{phieq}) with $\phi(\eta_1)=1$ and $\phi'(\eta_1)\leq0$ satisfies
\begin{equation}
\label{claim:L}
\phi'(\eta_2) < -1 \qquad \text{for all } \eta_2>\eta_1 \text{ with } \phi(\eta_2)\in[0,1/2].
\end{equation}
\end{lemma}

\begin{coro}\label{lemetabar}
For every
$\eta_2\geq\bar\eta,$
we have
\begin{equation}
\label{ineq:phi_leq_1}
\phi(\xi;\eta_2,\epsilon,-\bar \epsilon)< 1\quad
\mbox{for}\quad \xi\in[\bar \eta,\eta_2]
\end{equation}
for all $0\leq \epsilon\leq 1/2,$ and $0<\bar \epsilon\leq 1.$
\end{coro}

{\bf Proof of Corollary~\ref{lemetabar}.}
A violation of~(\ref{ineq:phi_leq_1}) implies the existence of
$\eta_1\in[\bar\eta,\eta_2)$ with $\phi(\eta_1)=1$ and $\phi'(\eta_1)\leq0$;
then~(\ref{claim:L}) contradicts the
condition
$\phi'(\eta_2;\eta_2,\epsilon,-\bar\epsilon) = -\bar\epsilon \geq -1$.\eop

{\bf Proof of Lemma~\ref{lemetabar2}.}
We fix
$\eta_1\gg1$ and define the
variable $y=\xi/\eta_1$. Changing variables,  equation (\ref{phieq})
transforms into
\begin{eqnarray}\label{phi:1}
0=\frac 1{\eta_1^2}\left(\ddot\phi+\frac
{d+1}{y}\dot\phi\right)-\frac
y{2}\dot\phi+\frac 1{d}y\phi\dot\phi+\phi^2-\phi,\quad \mbox{for}\quad
y>1\\
\dot\phi(1)=-D\eta_1, \quad \phi(1)=1,\label{phi:2}
\end{eqnarray}
where $\dot{} =\frac d{dy}.$ Define
\[
y_0 = \sup\{y>1: \ddot\phi(y) < 0 \mbox{ and } \phi(y)>0\}
\]
and note that $y_0>1$ if $\eta_1$ is large.
On $[1,y_0]$, $\dot\phi\leq -D \eta_1$; therefore $y_0\leq 1+ 1/(D\eta_1)$,
and consequently, on $[1,y_0]$,
\[
\fl
-\frac1d (y-1)\phi\dot\phi \leq \frac1{dD\eta_1} \mod{\dot\phi}
\leq \frac14\left(\frac12-\frac1d\right) \mod{\dot\phi}
\qquad \mbox{if}\qquad \eta_1\geq \frac{8d}{d-2}\, \frac1{dD}.
\]
Similarly,
\[
-\frac{d+1}{y\eta_1^2}\dot\phi \leq \frac14\left(\frac12-\frac1d\right)
\mod{\dot\phi}
\qquad \mbox{if}\qquad \eta_1^2\geq \frac{8d}{d-2}\, (d+1).
\]
Therefore,
\[
\frac1{\eta_1^2}\ddot\phi\leq \frac12\left(\frac12-\frac1d\phi\right)
\dot\phi - \phi^2 + \phi
\qquad \mbox{on }[1,y_0],
\]
for large $\eta_1$. Estimating $\mod{\dot\phi}$ by $D\eta_1$, we
find
\[
\frac1{\eta_1^2}\ddot\phi\leq -\frac12\left(\frac12-\frac1d\right) D\eta_1
+\frac14
\qquad \mbox{on }[1,y_0],
\]
and since the right-hand side of this expression is negative for large
$\eta_1$ it follows
that $\ddot\phi<  0$ on $[1,y_0]$; therefore $y_0$ may be redefined as
\[
y_0 = \sup\{y>1: \phi(y)>0\}.
\]
It follows that on $[1,y_0]$,
\begin{eqnarray*}
\dot\phi(y) &\leq&  -D\eta_1  - \eta_1^2 \frac d4 \left[\left(\frac12
-\frac1d \phi\right)^2
    - \left(\frac{d-2}{2d}\right)^2\right] + \int_1^y(\phi-\phi^2) \\
&\leq& -D\eta_1  - \eta_1^2 \frac d4 \left[\left(\frac12 -\frac1d
\phi\right)^2
    - \left(\frac{d-2}{2d}\right)^2\right] + \frac1{4D\eta_1}.
\end{eqnarray*}
When $0\leq \phi(y)\leq1/2$, this expression is bounded from above by
$-\eta_1^2/64$ for
large $\eta_1$. In terms of the original variable $\xi$ we obtain
$\phi'(\xi)\leq -\eta_1/64$, thus proving the lemma.\eop

\subsection{Definition of $\rho$ in $\Rset$}

The general
idea is to use $\eta_0=\bar \eta$ as a reference  point. In
this way, owing to Corollary~\ref{lemetabar}, we can obtain the
required estimates for $\rho.$ It can happen, however, that the
function $\phi(\xi,;\eta,v,w)$ is not defined at $\xi=\bar\eta.$ In
such a situation, to define $\rho,$ we introduce functions
representing the intersection of $\phi(\cdot;\eta,v,w)$ with the
lines $\phi=0$ for $\eta<\bar\eta$ and  $\phi=1$ for
$\eta>\bar\eta.$ Thus it is useful to define the following subsets of $\Rset:$
\begin{description}
\item $\Rset_1=\{(\eta,v,w)\in \Rset\colon\
\mbox{$\phi(\xi;\eta,v,w)$ satisfies (i) in Lemma
\ref{claphi}}\}$; \item $\Rset_2=\{(\eta,v,w)\in \Rset\colon\
\mbox{$\phi(\xi;\eta,v,w)$ satisfies (ii) in Lemma
\ref{claphi}}\}$, with  \begin{description}
\item $\Rset_{2a}=\Rset_2\cap\{\eta\leq \bar\eta\}$ and
$\Rset_{2b}=\Rset_2\cap\{\eta> \bar\eta\};$
\end{description}
\item $\Rset_3=\{(\eta,v,w)\in \Rset\colon\
\mbox{$\phi(\xi;\eta,v,w)$ satisfies (iii) in Lemma
\ref{claphi}}\}$.
\end{description}
We treat the cases in turn.

{\bf Case $\Rset_3$.} Fix a point  $(\eta,v,w)\in \Rset_3.$ We
choose $\eta_0=\bar \eta,$ $v_0=\phi(\bar \eta;\eta,v,w)$ and
$w_0=\phi'(\bar \eta;\eta,v,w).$ Note that
this choice is well defined:
since $\phi(\xi)\xi^2\to C>0$ as $\xi\to
\infty,$ there exists $\eta_\epsilon>\bar \eta$ such that
$\phi(\eta_\epsilon;\eta,v,w)=\epsilon<1/2,$ and
$-\phi'(\eta_\epsilon;\eta,v,w)=\bar\epsilon<1,$ with  $\bar
\epsilon\sim\frac{2\epsilon}{\eta_\epsilon}.$ Then Corollary~\ref{lemetabar} 
implies that the solution
$\phi(\cdot;\eta,v,w)$ can be continued to $\bar\eta$,
even
if $\bar \eta <\eta$. Setting
$\rho(\eta_0,v_0,w_0)=\eta_0^{d+1}e^{-\eta_0^2/4},$ we find
(\textit{cf.}~(\ref{def:rho}))
\begin{equation}
\fl\quad\rho(\eta,v,w)
\label{rhoetabar}
=\eta^{d+1}e^{-\eta^2/4}\exp\left\{
\frac1d \int\limits_{\bar\eta}^\eta \xi
\phi(\xi;\bar\eta,\phi(\bar \eta;\eta ,v,w),\phi'(\bar
\eta;\eta,v,w))\, d\xi\right\}.
\end{equation}
The choice of $\eta_0=\bar \eta$ also allows us to estimate the value
of $\phi$ for $\xi>\bar \eta,$ which in turn permits us to control
$\rho$ for large $\eta$, since the bound
$\phi(\xi)\leq 1$ for $\xi>\bar \eta$
implies an exponential decay for $\rho$ as $\eta\to
\infty.$

{\bf Case $\Rset_1$.} Points in  $\Rset_1$ are of the form
$(\eta,1,0)$ and $(\eta,0,0).$ We again choose
$\eta_0 = \bar\eta$; substituting  $\phi\equiv 1$ and
$\phi\equiv 0$ into formula (\ref{rhoetabar}) gives
\begin{eqnarray}\label{rho1.0}
\rho(\eta,1,0)=\eta^{d+1}e^{-\frac{(d-2)\eta^2}{4d}}e^{-\frac{\bar\eta^2}{2d}},\quad\mbox{and}\quad
\rho(\eta,0,0)=\eta^{d+1}e^{-\frac{\eta^2}{4}}.
\end{eqnarray}

{\bf Case $\Rset_{2a}.$}  Fix a point $(\eta,v,w)\in \Rset_{2a}.$
Let $\eta^*$ be given by Lemma \ref{claphi} and  define the
function $L_0\colon \Rset_{2a}\to \R^+$ such that
$L_0(\eta,v,w)=\min\{\eta^*,\bar \eta\}.$ Note that the function
$L_0$ is continuous and equals either the point $\eta^*$ where
$\phi(\eta^*;\eta,v,w)$ vanishes or $\bar \eta$ if
$\phi(\bar \eta;\eta,v,w)\geq 0.$
To find  $\rho,$
we choose $(\eta_0,v_0,w_0)=(\eta^*,0,\phi'(\eta^*,\eta,v,w)),$
and set $\rho(\eta_0,v_0,w_0)=\eta_0^{d+1}e^{-\eta_0^2/4}.$
This gives
\begin{equation}
\rho(\eta,v,w)=\eta^{d+1}\exp\{-\eta^2/4+I_0\},
\label{rhoL}
\end{equation}
where
$$
I_0=\int\limits_{
L_0(\eta,v,w)}^\eta \frac 1{d}\xi \phi\left(\xi;
L_0(\eta,v,w),\phi( L_0(\eta,v,w);\eta,v,w),\phi'( L_0(\eta,v,w);\eta,v,w)\right)\,
d\xi.
$$
{\bf Case $\Rset_{2b}.$} Here it is convenient to define for any
$(\eta,v,w)\in R_{2b}$ the function $L_1\colon \Rset_{2b}\to
\R^+,$ by
\begin{equation}
\label{def:L1}
\fl \quad L_1(\eta,v,w)= \cases{ 
 \max\{\bar\eta,\max\{\xi \in (0,\eta) \mid
\phi(\xi;\eta,v,w)\geq 1\}\}\quad&if $v<1$,\\ \\
\min\{\xi \in (\eta,\infty) \mid \phi(\xi;\eta,v,w)\leq
1\}\quad&if $v\geq 1$. }
\end{equation}
The function $L_1$ is well defined
for $v<1$ since if
$\phi(\tilde\xi;\eta,v,w)=0$ for some  $\tilde \xi\in (\bar
\eta,\eta)$ then $\phi<1$ in $(\bar \eta,\eta)$ by
Corollary~\ref{lemetabar} and $\phi$
has to attain
a local maximum in  $(\bar \eta,\eta)$, which is a contradiction with equation~(\ref{phieq}).
For $v\geq 1$, $L_1$ is well-defined by Lemma \ref{claphi}.

Note that $\phi(L_1(\eta,v,w);\eta,v,w)\leq 1.$ The function
$L_1$ is
continuous and equals either $\eta_*$ where
$\phi(\eta_*;\eta,v,w)=1$ or $\bar \eta$ if $\phi(\bar
\eta;\eta,v,w)\in (0,1)$.

Now fix a point $(\eta,v,w)\in \Rset_{2b},$  choose $\eta_0=
L_1(\eta,v,w)$ and set
$\rho(\eta_0,v_0,w_0)=\eta_0^{d+1}e^{-(d-2)\eta_0^2/4d}e^{-\bar
\eta^2/2d}.$ Using (\ref{def:rho}),  we find
\begin{equation}\label{rhoL1a}
\rho(\eta,v,w)=\eta^{d+1}\exp\{-\eta^2/4+\eta_0^2/2d-{\bar
\eta}^2/2d+I_1\}
\end{equation}
where
$$
I_1=\int\limits_{L_1(\eta,v,w)}^\eta \frac 1{d}\xi \phi\left(\xi,
L_1(\eta,v,w),\phi( L_1(\eta,v,w);\eta,v,w),\phi'( L_1(\eta,v,w),\eta,v,w)\right)\,
d\xi, 
$$
and 
\begin{equation}
\label{rhoL1}
\rho(\eta,v,w)=\eta^{d+1}\exp\{{-(d-2)\eta^2/4d}-{\bar
\eta}^2/2d+{I'_1}\}
\end{equation}
where
$$
I'_1=
\int\limits_{L_1(\eta,v,w)}^\eta \frac 1{d}\xi [\phi\left(\xi,
L_1(\eta,v,w),\phi( L_1(\eta,v,w);\eta,v,w),\phi'( L_1(\eta,v,w),\eta,v,w)\right)-1]\,
d\xi.
$$

\subsection{Properties of $\rho$ and $\Phi$}
In the previous section, we have found a solution $\rho$ of
(\ref{f.ord.rho}). Here we show that this  solution, together with
the function $\Phi$ given by~(\ref{pPhi2}),
satisfies the properties required for the proof of Theorem
\ref{converg}. We start by  stating a result which provides a lower
bound for $\rho$ in $\Rset_{2b}.$
\begin{lemma}\label{Gfunction} Let  $M$ and $\bar M$ be the
constants in estimates
(\ref{boundB}) and (\ref{boundBeta}), and let $L_1$ be defined as in~(\ref{def:L1}). Then there
exists a large constant  $\bar \eta_0$ such that the function
$G\colon[\bar \eta_0,\infty)\to \R^+$ given by
\[
\fl\qquad G(\eta)=\max\{ L_1(\eta,a,-b)\mid 1\leq a\leq M\:\mbox{and}\: 0\leq b\leq \bar
M\:\}\quad \mbox{for}\quad \eta\geq \bar \eta_0,
\]
satisfies $G(\eta)\leq C\eta$ for some constant $C=C(M)>0.$
\end{lemma}
{\bf Proof.} We take $\bar \eta_0$ large and we fix $\eta\geq\bar \eta_0.$
Using the  continuity of $L_1,$
we have that $G(\eta)= L_1(\eta,\bar a,-\bar b)$  for some $\bar a\in [1,M],$ and $\bar b\in [0,\bar M].$ Now we define the variable $y=\xi/\eta\geq 1$; the result is proved if
we show that $\sup \{y\geq 1: \phi(y)> 1\} \leq C(M)$.

As in the proof of Lemma~\ref{lemetabar2}, equation (\ref{phieq})
transforms into
\begin{eqnarray}\label{phidot2a}
\fl\quad 0=\frac
1{\eta^2}\left(\ddot\phi+\frac{d+1}{y}\dot\phi\right)-\frac
y{2}\dot\phi+\frac 1{d}y\phi\dot\phi+\phi^2-\phi \quad
\mbox{for}\quad y>1,\\ \quad \dot\phi(1)=-\bar b\eta
\:\:\mbox{and}\:\: \phi(1)=\bar a.
\end{eqnarray}
Note that for $\phi>1$ we have $\dot\phi(y)<0$ for all $y>1,$ since
$\dot \phi(\bar y)=0$ implies that $\bar y$ can only be a maximum, a contradiction with
equation~(\ref{phieq}).

We prove the claim in two steps. In the first step we consider the case $\bar a> d/2 - \delta>1$,
where $\delta = (d+1)d/\eta^2$.
Define $y_1 = \sup\{y>1: \phi(y)>d/2-\delta\}$.
We write (\ref{phidot2a}) as 
\begin{equation}\label{inty2a}
\frac 1{\eta^2}\ddot\phi= -yA_2(y)\dot\phi-A_1(y)
\quad \mbox{for}\quad y>1,
\end{equation}
where 
\[
A_1(y)=\phi^2-\phi\quad \mbox{and}\quad A_2(y)=\left(\frac 1{d}\phi-\frac 12+\frac{d+1}{y^2\eta^2}\right).
\]
Since $\phi(\cdot)\in [d/2-\delta,\bar a]$ on $[1,y_1]$, 
$A_2$ is non-negative and bounded by  $\bar A_2:=\bar a/d.$  
The function $A_1$ is positive and bounded from below:
\[
A_1(y) \geq \underline{A_1} := \left(\frac d2-{\delta} \right)^2 -\left(\frac d2-\delta \right) > 0.
\]
Integrating equation (\ref{inty2a}), we have
\begin{equation}\label{inty2b}
\fl\quad
\dot \phi(y)=-\bar b\eta \,e^{-\eta^2\int\limits_1^ytA_2(t)\,dt} 
  -\eta^2\int\limits_1^yA_1(s)e^{-\eta^2\int\limits_s^ytA_2(t)\,dt}\,ds \quad \mbox{for}\quad y>1.
\end{equation}
We observe that 
\begin{equation}\label{inty2c}
\eta^2\int\limits_1^yA_1(s)e^{-\eta^2\int\limits_s^ytA_2(t)\,dt}\,ds\geq
\underline{A_1}f(y;\eta) \quad \mbox{for}\quad 1\leq y\leq y_1
\end{equation}
where 
$
f(y;\eta)=\eta^2\int\limits_1^y e^{-\eta^2\bar A_2(y^2-s^2)/2}\,ds
$
is a positive bounded function satisfying $yf(y;\eta)\to 1/\bar A_2$ as $y\to \infty$
(the latter claim follows from considering the integrand close to $s=y$), and 
more precisely,
\[
yf(y;\eta) \geq \frac 1{2\bar A_2} 
\qquad \text{for $y\geq 2$ and for sufficiently large $\eta$}.
\]
Therefore the primitive function
\[
F(y; \eta) = \int_1^y f(s;\eta)\, ds
\]
satisfies
\begin{equation}
\label{bound:F}
F(y;\eta) \geq \frac 1{2\bar A_2}(\log y -\log 2).
\end{equation}
Integrating (\ref{inty2b}) on $[1,y_1]$ and using (\ref{inty2c}),  we obtain  
\[
\phi(y_1)\leq \bar a-\bar b\eta(y_1-1)- \underline{A_1} F(y_1;\eta).
\]
To obtain a bound on $y_1$, we use $\phi(y_1) = d/2 -\delta$ and conclude
\[
\underline{A_1}F(y_1;\eta) \leq \bar a \leq M,
\]
from which it follows that
$
y_1 \leq C(M)
$
by~(\ref{bound:F}).

\medskip

For the second step, we replace $\eta$ by $y_1\eta$ in the rescaling above,
by which we can assume that we are in the same situation: $\phi(1) = \bar a$,
$\dot\phi(1) = \bar b \eta$, but this time $1\leq \bar a\leq d/2-\delta$. 

Similarly define $y_2 = \sup\{y\geq1: \phi(y)> 1\}$. Since $1\leq \phi(\cdot)\leq d/2-\delta$
on $[1,y_2]$, the function $A_2(\cdot)$ in~(\ref{inty2a}) is negative,
so that $\phi$ satisfies the differential inequality
\begin{equation}
\label{comp:psi}
\frac1{\eta^2}\ddot\phi  \leq -\phi^2 +\phi < -2(\phi-1).
\end{equation}
Let the function $\psi$ solve
\[
\frac1{\eta^2} \ddot \psi = -2(\psi-1), \qquad \text{with } \psi(1) = \bar a 
\text{ and }\dot\psi(1) = 0.
\]
The solution of this equation is $\psi(y) = 1+ (\bar a-1) \cos(\eta\sqrt 2(y-1))$,
and note that $\psi(\tilde y_2)=1$ for $\tilde y_2 := \pi/(2\eta\surd 2)$.
From~(\ref{comp:psi}), $\phi(1+) < \psi(1+)$; if $\phi(y) = \psi(y)$ for some $y\in(1,\tilde y_2)$,
then by the comparison principle (which the operator $u\mapsto \ddot u/\eta^2 + 2(u-1)$
satisfies on intervals of length less than $\tilde y_2$) we find $\phi\geq\psi$ on
the interval $[1,y]$, a contradiction
with the previous remark.

In conclusion we find that $y_2\leq \tilde y_2$, thus proving the lemma.
\eop 

\bigskip

We now derive estimates for $\rho$ and $\Phi$ in $\tilde \Rset$
and  $\Rset.$
\begin{lemma}\label{prho}
The function $\rho$ is continuous in $\Rset\setminus
\{\eta=\bar\eta,\: v>1\}$; for  $(\eta,v,w)\in \Rset,$ one
finds
\begin{eqnarray}\label{prhoe.a}
\rho(\eta,v,w)\leq \eta^{d+1} e^{-(d-2)\eta^2/4d}.
\end{eqnarray}
In addition, if $(\eta,v,w)\in \tilde \Rset,$ then
\begin{eqnarray}\label{prhoe.b}
\rho(\eta,v,w)\geq \frac 1{C_0}\eta^{d+1} e^{-C_0\eta^2}
\end{eqnarray}
for some constant $C_0=C_0(M)>0.$
\end{lemma}
{\bf Proof.} We start by proving (\ref{prhoe.a})--(\ref{prhoe.b}).
Let $\tilde \Rset_i=\tilde \Rset\cap \Rset_i$ for $i=1,2,3.$ 
If $(\eta,v,w)\in \Rset_1$ then the estimates
(\ref{prhoe.a})--(\ref{prhoe.b}) follow by definition. 
If $(\eta,v,w)\in \Rset_{2a},$ then as  $\phi>0$ on  $(\eta,L_0(\eta,v,w))$
the integral in (\ref{rhoL}) is negative. This
gives $$\rho(\eta,v,w)\leq  \eta^{d+1} e^{-\eta^2/4}\quad \mbox{for}\quad(\eta,v,w)\in \Rset_{2a}.
$$ 
Now for $(\eta,v,w)\in \tilde\Rset_{2a},$ we have that $[\eta, L_0(\eta,v,w)]\subset[0,\bar \eta],$ $v\in[0,M]$ and $w\in [-\bar M,0].$ Then the continuity of $\phi$ on $[\eta, L_0(\eta,v,w)]$ implies that 
\begin{equation*}
\fl|\phi\left(\cdot;
L_0(\eta,v,w),\phi( L_0(\eta,v,w);\eta,v,w),\phi'( L_0(\eta,v,w);\eta,v,w)\right)|_{C^0([\eta, L_0(\eta,v,w)])}\leq \bar C_0
\end{equation*}
where $\bar C_0=\bar C_0(M,\bar M,\bar \eta).$  Using this bound to estimate $I_0$ in (\ref{rhoL}), we find 
$$
C(M)\eta^{d+1}
e^{-\eta^2/4}\leq \rho(\eta,v,w)\quad\mbox{for}\quad (\eta,v,w)\in \tilde\Rset_{2a}
$$ 
with $C(M)<1,$ since we have integrated backwards.

For any $(\eta,v,w)\in \Rset_{2b},$ we use (\ref{rhoL1}) and find
the upper bound $$\rho(\eta,v,w)\leq \eta^{d+1} e^{-(d-2)\eta^2/4d}
e^{-\bar\eta^2/2d}.$$ This estimate follows from the negative sign  of the integral $I'_1$ in (\ref{rhoL1}). In fact, for  $\xi\in (\eta,L_1(\eta,v,w)),$ we have 
$$
\cases{\phi\left(\xi;
L_1(\eta,v,w),\phi( L_1(\eta,v,w);\eta,v,w),\phi'( L_1(\eta,v,w);\eta,v,w)\right)-1>0\quad \mbox{if $v\geq 1,$} \cr  
\phi\left(\xi;
L_1(\eta,v,w),\phi( L_1(\eta,v,w);\eta,v,w),\phi'( L_1(\eta,v,w);\eta,v,w)\right)-1<0\quad\mbox{if $v<1.$}}
$$
Next for $(\eta,v,w)\in \tilde\Rset_{2b},$ we find 
$$
\cases{ 
\rho(\eta,v,w)\geq  \eta^{d+1} e^{-\eta^2/4}& $\quad
\mbox{for}\quad v\leq 1,$   \cr \rho(\eta,v,w)\geq  \eta^{d+1}
e^{-(d-2)\eta^2/4d} e^{-\bar\eta^2/2d}e^{-\bar C(M)\eta^2}&$\quad
\mbox{for}\quad v>1,$
}
$$
where $\bar C(M)>0.$ The estimate when $v\leq 1$ follows directly from (\ref{rhoL1a}). To obtain the estimate for $\rho$ when $v>1$, we use (\ref{rhoL1}). In fact, noting that $\phi$ is   
non increasing in $[\eta, L_1(\eta,v,w)],$ we have that 
\begin{equation*}
\fl|\phi\left(\cdot;
L_1(\eta,v,w),\phi( L_1(\eta,v,w);\eta,v,w),\phi'( L_1(\eta,v,w);\eta,v,w)\right)|_{C^0([\eta, L_1(\eta,v,w)])}\leq M.
\end{equation*}
Using this bound together with the estimate $L_1(\eta,v,w)\leq C(M)\eta$ (see  Lemma \ref{Gfunction}), we find that $I'_1$ in (\ref{rhoL1}) satisfies 
$
-I'_1\leq \bar C(M)\eta^2,
$
which gives the derired estimate.

To prove (\ref{prhoe.a}) for $(\eta,v,w)\in \Rset_3,$ we examine
two cases, if $\eta\leq\bar\eta$ then the estimate for
$\Rset_{2a}$ holds and for $\eta>\bar\eta$ the estimate for
$\Rset_{2b}$ holds. Finally, to obtain (\ref{prhoe.b}) for
$(\eta,v,w)\in \tilde\Rset_3,$ we also check two cases, if
$\eta\leq\bar\eta$ then the estimate for $\tilde\Rset_{2a}$ holds
and for $\eta>\bar\eta$ the estimate for $\tilde\Rset_{2b}$ with
$v\leq 1$ holds.

{\bf Claim.} $\rho$ is continuous in $\Rset\setminus
\{\eta=\bar\eta,\: v>1\}.$

Before we prove this, note that $\Rset_2$ is an open set  and
$\Rset_1$ and $\Rset_3$ are closed.

We first see that $\rho$ is continuous within  $\Rset_{2a}$ and
$\Rset_{2b}$, by  continuity of $L_0$ and $L_1.$ For the elements
in $\Rset_1,$ the definition of $\rho$ is as for $\Rset_2,$
therefore there is continuity of $\rho$ between $\Rset_2$ and
$\Rset_1.$

The delicate part is to proof continuity between $\Rset_3$ and
$\Rset_2$. Taking a sequence $(\eta_n,v_n,w_n)\in \Rset_2,$ we
associate a solution $\phi_n(\cdot,\eta_n,v_n,w_n).$ Suppose that
$(\eta_n,v_n,w_n)\to (\eta,v,w)\in \Rset_3.$ Now if
$\phi(\cdot,\eta,v,w)$ is the solution  of (\ref{phieq}) then
$\phi_n\to\phi$ in compact subsets of $\R^+.$ Therefore by Corollary~\ref{lemetabar}, 
for  $n\geq n_0\in \N,$ we find
$\phi_n(\bar\eta)\in(0,1).$
Then $(\eta_n,v_n,w_n)\in \Rset_{2}$ for $n\geq n_0,$ have the
same definition of $\rho$ as  for $(\eta,v,w)\in \Rset_3.$ Finally
if $v\leq 1$ and $\eta=\bar \eta$, then $\rho$ is continuous. If
$\eta$ close enough to $\bar \eta$ then we have that  $\eta_0=\bar
\eta$. So the computation of $\rho$ uses  the same formula,
independent of the subset of $\Rset$ to which $(\eta,v,w)$
belongs. \eop

For $\Phi$ we deduce the following lemma, which implies (\ref{Phiest}).
\begin{lemma}\label{pPhi}
The function $\Phi$ is continuous in $\Rset\setminus\{\eta=\bar
\eta,\: v>1\}$ and if $(\eta,v,w)\in R$, then
$$
\Phi(\eta,v,w)\leq \left\{w^2 +\frac{v^2}{2}\right\} \eta^{d+1} e^{-(d-2)\eta^2/4d}
$$
and
$$
\Phi(\eta,v,w)\geq -\left\{\frac{v^3}{3}-\frac{v^2}{2}\right\}\eta^{d+1} e^{-(d-2)\eta^2/4d}.
$$
\end{lemma}
{\bf Proof.} Follows directly from the definition  (\ref{pPhi2})
of $\Phi$ and uses the upper bound  (\ref{prhoe.a}) of $\rho$.
\eop
\subsection{Regularizing argument}\label{regularPhi}
In the beginning of this appendix, we formally constructed  a
Lyapunov functional $E(\tau)$ with $\Phi$ and $\rho$ satisfying
(\ref{lyapbound}). In the previous section, we obtained a solution
$\rho$ of  (\ref{f.ord.rho}) and $\Phi$ given
 by (\ref{Phicond}). Moreover these functions satisfy the properties
 found in  Lemmas \ref{prho} and \ref{pPhi}. From these
results
we do not obtain enough regularity to derive (\ref{lyapbound}). To
do this, we introduce a  introduce a regularization of $\Phi$
using standard mollifiers and translation function to avoid the
singularity of $f$ at $\eta=0.$ See the details of the proof in
\cite[p. 102]{g}.
\section{Appendix: Linear stability of blow-up profiles}\label{stability}

In this appendix, we study the linear stability of the blow-up
profiles $\varphi_1$ and $\varphi^*,$ see (\ref{phis1}).

Let $B$ be a solution of (\ref{equaB1})--(\ref{equaB3}) and let
$\varphi$ be a solution of (\ref{staW.n}). The idea is to study
the linearized equation for the difference
$\Phi(\eta,\tau):=B(\eta,\tau)-\varphi(\eta),$ i.e.
\begin{equation}\label{linearpart}
\fl\qquad\Phi_\tau=\Phi_{\eta\eta}+\frac{d+1}{\eta}\Phi_\eta+\left(\frac
1{d}\varphi-\frac 1{2}\right)\eta\Phi_\eta+\left(\frac
1{d}\eta\varphi_\eta+2\varphi-1\right)\Phi.
\end{equation}
Here, we have implicitly assumed that sufficiently close to
blow-up only the linear terms play a role in describe the
singularity formation.

For the stability analysis, let $\lambda>0$ and consider
a solution  of (\ref{linearpart}) of the form $\psi_{\lambda}(\eta)e^{\lambda\tau}.$ By (\ref{linearpart}), $\langle\psi_{\lambda}(\eta),\lambda\rangle$ satisfies
\begin{equation}\label{linearpart2}
\fl\qquad
(\psi_{\lambda})_{\eta\eta}+\frac{d+1}{\eta}(\psi_{\lambda})_\eta+\left(\frac
1{d}\varphi-\frac 1{2}\right)\eta(\psi_\lambda)_\eta+\left(\frac
1{d}\eta\varphi_\eta+2\varphi-1-\lambda\right)\psi_\lambda.
\end{equation}
For the analysis of boundary conditions we consider first
$\varphi=\varphi_1.$
We note that at $\eta=0$ we have either $\psi_{\lambda}\sim 1$ or $\psi_{\lambda}\sim 1/\eta^d$. To have $\psi_{\lambda}$ bounded near 0, we impose
\begin{equation}\label{bcl1}
(\psi_{\lambda})_\eta(\eta)\to 0, \quad \mbox{as} \quad \eta\to 0.
\end{equation}
For large $\eta,$ we can either have $\psi_{\lambda}\sim \eta^{-(2\lambda+3)}e^{\eta^2/4}$ or $\psi_{\lambda}\sim \eta^{2\lambda-2}$. We see that both behaviours diverge with $\eta$, however the second asymptotic is bounded in terms of $r$ and $t$ as $t\to T.$ Therefore to have a  polinomial behaviour at infinity, we  prescibe
\begin{equation}\label{bcl2}
\psi_{\lambda}(\eta)e^{-\eta}\to 0, \quad \mbox{as} \quad \eta\to \infty.
\end{equation}
Now solving  equation (\ref{linearpart2}) together with (\ref{bcl1}) and (\ref{bcl2}),  we find a sequence
solutions of (\ref{linearpart}) given by
$\{e^{\lambda_n\tau}\psi_{n}(\eta)\}_{n\in \N\cup\{0\}},$ with
$\lambda_0>\lambda_1>\ldots,$ where $\psi_{n}:=\psi_{\lambda_n}.$
If the blow-up time $T>0$ is chosen correctly in the
definition of $\eta$ and $\tau,$ we can eliminate, see
\cite{bcksv}, the first mode ($n=0$) corresponding to change of blow-up  and write
\[
B(\eta,\tau)=\varphi(\eta)+\psi_1(\eta)e^{\lambda_1\tau}+O(e^{\lambda_2\tau}).
\]
Therefore from the sign of $\lambda_1$ we obtain the linear stability
of $\varphi$.


In \cite{bcksv}, Brenner et al. proved, using (\ref{linearpart}),
the following stability result  for various blow-up profiles.
\begin{theorem}
Every solution $\varphi$ of (\ref{staW2.a}) satisfying $\eta
\varphi_\eta/\varphi\to 2$ as $\eta\to \infty$ has an unstable
mode corresponding to changing the blow-up time. Also, a blow-up
profile with $k$ intersections with the singular solution $\varphi_S$
has at least $k-1$ additional unstable modes.
\end{theorem}
In addition, the authors in \cite{bcksv} found numerically that $\lambda_1<0$ when
$\varphi=\varphi_1$ and $d>2.$
In particular, they computed $\lambda_1=-0.272\ldots$ for $d=3.$  This implies that
$\varphi_1$ is linearly stable for $d>2.$

For $\varphi=\varphi^*,$ we can proceed as above and solve the
eigenvalue problem for (\ref{linearpart}).
Considering (\ref{linearpart2}) with  $\varphi=\varphi^*,$  we find that $\langle \psi_\lambda,\lambda\rangle$ satisfies
\begin{equation}\label{eq:phiast}
(\psi_\lambda)_{\eta\eta}+\left(\frac {d+1}{\eta}-\frac
{d-2}{2d}\eta\right)(\psi_\lambda)_\eta+(1-\lambda)\psi_\lambda=0,
\end{equation}
with (\ref{bcl1}) and (\ref{bcl2}). These boundary conditions are chosen by the same arguments for $\varphi=\varphi_1;$ however in the current case
we have either $\psi_\lambda\sim
\eta^{\frac{2d}{d-2}(\lambda-1)-d-2}e^{\frac{(d-2)}{4d}\eta^2}$ or
$\psi_\lambda\sim \eta^{\frac{2d}{d-2}(1-\lambda)}$ as $\eta\to\infty.$  Note that by changing $\eta$ by $(-\eta)$ the equation remains invariant, so only solutions consisting on even
powers are allowed. Then we construct a sequence of solutions of the form
\[
\psi_n(\eta)=\sum\limits_{i=0}^{n} A_i \eta^{2i}\quad \mbox{for any}\quad n=0,1,2,3\ldots,
\]
where the coefficients are given by $A_{i}(2i(2i-1)+(d+1)2i)=A_{i-1}(1-\lambda-2i(d-2)/2d)$ for $i=1,2,\ldots$ and $A_0$ an  arbitrary constant. This means that when $(1-\lambda-2(n+1)(d-2)/2d)=0,$ we find an explicit polinomial solution of degree $2n$, where $\lambda$ is given by
\begin{equation}\label{lambdan}
\lambda_n=\frac{d-n(d-2)}{d}.
\end{equation}
Consequently, we have obtained an explicit sequence of solution $\{\langle \psi_n,\lambda_n\rangle\}_{n\in \N\cup\{0\}}$ for the eigenvalue problem (\ref{eq:phiast}). The eigenvalue $\lambda_0=1$ corresponds to the unstable mode of
change of blow-up time and since  $\lambda_1>0$ for all $d>2,$ by (\ref{lambdan}), this means that $\varphi^*$ is linearly
unstable.

\end{document}